\documentclass[11pt,reqno]{amsart}

\usepackage{amscd,amsmath,amssymb,amsfonts}
\usepackage{stmaryrd}
\usepackage{bbold}
\usepackage{times}
\usepackage[all]{xy}
\usepackage{mathrsfs}

\usepackage{amsthm}

\theoremstyle{plain}
\newtheorem{thm}{Theorem}
\newtheorem{lem}[thm]{Lemma}
\newtheorem{cor}[thm]{Corollary}
\newtheorem{prop}[thm]{Proposition}

\newtheorem*{prop*}{Proposition}

\theoremstyle{definition}

\theoremstyle{definition}

\newtheorem{example}[thm]{Example}
\newtheorem*{example*}{Example}
\newtheorem*{ack}{Acknowledgement}

\numberwithin{thm}{section}
\numberwithin{equation}{section}

\newcommand{\si}{\sigma}
\newcommand{\e}{\varepsilon}

\newcommand{\onto}{\twoheadrightarrow}
\newcommand{\mono}{\rightarrowtail}
\newcommand{\iso}{\overset{\sim}{\tto}}

\newcommand{\Hom}{\operatorname{Hom}}
\newcommand{\Aut}{\operatorname{Aut}}
\newcommand{\im}{\operatorname{Im}}

\DeclareMathOperator{\Id}{Id}
\DeclareMathOperator{\ch}{char}
\DeclareMathOperator{\Tor}{Tor}

\DeclareMathOperator{\Ext}{Ext}

\newcommand{\Mat}{\operatorname{Mat}}

\newcommand{\sA}{{\mathcal A}}
\newcommand{\sB}{{\mathcal B}}
\newcommand{\sC}{{\mathcal C}}

\newcommand{\sO}{{\mathcal O}}

\newcommand{\sR}{{\mathcal R}}
\newcommand{\sS}{{\mathcal S}}

\newcommand{\sX}{{\mathcal X}}

\newcommand{\Z}{{\mathbb Z}}

\newcommand{\q}{\mathbf{q}}
\newcommand{\bg}{\mathbf{g}}
\newcommand{\bi}{\mathbf{i}}
\newcommand{\bj}{\mathbf{j}}
\newcommand{\T}{\mathsf{T}}
\newcommand{\K}{\mathsf{K}}
\newcommand{\ju}{\nu}
\newcommand{\Sym}{\mathsf{S}}
\newcommand{\Gr}{\mathsf{\Lambda}}
\newcommand{\YM}{{\mathcal Y\negmedspace \mathcal M}}

\newcommand{\Sets}{\text{\sf Sets}}

\newcommand{\Vect}{\text{\sf Vect}}
\newcommand{\HN}{\text{\sf H}_N\text{\sf Alg}^s_{\k}}
\newcommand{\comod}{\text{\sf Comod}^s}
\newcommand{\com}{\text{\sf comod}^s}

\newcommand{\ihom}{\operatorname{\underline{Hom}}}
\renewcommand{\hom}{\operatorname{\underline{hom}}}
\newcommand{\en}{\operatorname{\underline{end}}}
\newcommand{\eA}{\en \sA}

\newcommand{\tto}{\longrightarrow}
\newcommand{\into}{\hookrightarrow}
\newcommand{\dg}{\widehat}
\DeclareMathOperator{\Ker}{Ker}
\DeclareMathOperator{\sdim}{sdim}
\newcommand{\schi}{\chi^s}
\renewcommand{\k}{\mathbb{k}}
\newcommand{\f}{\mathbf{f}}
\newcommand{\Sy}{\mathfrak{S}}

\newcommand{\He}{\mathscr{H}}
\newcommand{\scrT}{\mathscr{T}}

\newcommand{\Red}{\operatorname{Red}}
\newcommand{\NRed}{\operatorname{NRed}}
\newcommand{\me}{\wedge}
\newcommand{\jo}{\vee}
\newcommand{\sign}{\operatorname{sgn}}

\def\End{\operatorname{End}}
\def\Endo{\text{\sffamily E}}
\def\inv{\text{\sffamily inv}}
\def\GL{\text{\sffamily GL}}

\def\ev{\text{\rm ev}}
\def\ber{\operatorname{ber}}
\def\tr{\text{\rm tr}}
\def\str{\text{\rm str}}

\newcommand{\hR}{{\mathscr R}}

\newdir{ >}{{}*!/-5pt/\dir{>}}

\begin{document}

\title[$N$-homogeneous superalgebras]{$N$-homogeneous superalgebras}

\author{Ph{\`u}ng H{\^o} Hai}
\address{Mathematik, University of Duisburg-Essen, Germany and
Institute of Mathematics, Hanoi, Vietnam}
\email{hai.phung@uni-duisburg-essen.de}
\thanks{PHH is supported by the DFG through a Heisenberg-Fellowship}

\author{Benoit Kriegk}
\address{LaMUSE, Facult{\'e} des Sciences et Techniques,
Universit{\'e} de Saint-Etienne, 23 rue du Docteur Paul Michelon,
42023 Saint-Etienne Cedex 2, France}
\email{benoit.kriegk@univ-st-etienne.fr}
\thanks{}

\author{Martin Lorenz}
\address{Department of Mathematics, Temple University,
    Philadelphia, PA 19122-6094, USA}
\email{lorenz@temple.edu}
\thanks{ML's research is supported in part by NSA Grants H98230-05-1-0025
and H98230-07-1-0008 and by Leverhulme Research Interchange Grant
F/00158/X}



\maketitle


\begin{abstract}
We develop the theory of $N$-homogeneous algebras in a super
setting, with particular emphasis on the Koszul property. To any
Hecke operator $\hR$ on a vector superspace, we associate certain
superalgebras $\Sym_{\hR,N}$ and $\Gr_{\hR,N}$ generalizing the
ordinary symmetric and Grassmann algebra, respectively. We prove
that these algebras are $N$-Koszul. For the special case where $\hR$
is the ordinary supersymmetry, we derive an $N$-generalized
super-version of MacMahon's classical ``master theorem".
\end{abstract}


\setcounter{tocdepth}{1}
\tableofcontents


\section*{Introduction}

\subsection{} \label{SS:background}

The theory of $N$-homogeneous algebras owes its existence primarily
to the concerns of noncommutative geometry. In fact, as has been
expounded by Manin in his landmark publications \cite{yM87},
\cite{yM88}, quadratic algebras (the case $N=2$) provide a
convenient framework for the investigation of quantum group actions
on noncommutative spaces. Moreover, certain Artin-Schelter regular
algebras \cite{mAwS87}, natural noncommutative analogs of ordinary
polynomial algebras, can be presented as associative algebras
defined by cubic relations ($N=3$). The latter algebras, as well as
many of the quadratic algebras studied by Manin, enjoy the
additional ``Koszul property" which will be of central importance in
the present article; it will be reviewed in detail in
\ref{SS:generalized} below.

Motivated by these examples and others, Berger \cite{rB01} initiated
the systematic investigation of $N$-homogeneous algebras for all $N
\ge 2$, introducing in particular a natural extension of the notion
of Koszul algebra from the familiar quadratic setting to general
$N$-homogeneous algebras. Article \cite{rB01} gives examples of
$N$-Koszul algebras for all $N \ge 2$; these are the so-called
$N$-symmetric algebras, the special case $N=2$ being the ordinary
symmetric (polynomial) algebra. Following the general outline of
Manin's lecture notes \cite{yM88} on the case of quadratic algebras,
Berger, Dubois-Violette and Wambst developed the categorical aspects
of $N$-homogeneous algebras in \cite{rBmDVmW}.

\subsection{} \label{SS:superization}

Current interest in $N$-homogeneous algebras is fueled in part by
the fact that they do occur naturally in mathematical physics and in
combinatorics. Indeed, Connes and Dubois-Violette \cite{aCmDV02},
\cite{aCmDV03} introduced a class of $3$-homogeneous algebras,
called Yang-Mills algebras, which are in fact $3$-Koszul. There are
two versions of Yang-Mills algebras: in the language of linear
superalgebra, the first kind has even (parity $\bar 0$) algebra
generators while the second kind is generated by odd (parity $\bar
1$) elements.

Combinatorics enters the picture via MacMahon's celebrated ``master
theorem" \cite{paMM60}, specifically the recent quantum
generalization of the master theorem due to Garoufalidis, L{\^e} and
Zeilberger \cite{GLZxx}. As has been pointed out by two of the
present authors in \cite{HL}, the yoga of (quadratic) Koszul
algebras leads to a rather effortless and conceptual proof of the
quantum master theorem based on the fact that a certain quadratic
algebra, known as quantum
affine space, is Koszul. 
Further quantum generalizations and super versions of the master
theorem have been obtained by several authors using a variety of
approaches; see Foata and Han \cite{FH06a}, \cite{FH06b},
\cite{FH06c}, Konvalinka and Pak \cite{mKiP}, Etingof and Pak
\cite{EP06}.

\subsection{} \label{SS:MT}

From an algebraic point of view, MacMahon's master theorem (MT) in
its various incarnations finds its most natural explanation by the
phenomenon of ``Koszul duality". Indeed, all versions of MT can be
expressed in the form that, for some algebra $\sB$, an equation
$\Sigma_1\cdot \Sigma_2 = 1$ holds for suitable power series
$\Sigma_1, \Sigma_2 \in \sB\llbracket t \rrbracket$. Here is a brief
outline how one can arrive at such an equation starting with a given
$N$-Koszul algebra $\sA$. Associated with $\sA$, there is a graded
complex, $\K(\sA)$, which is exact in positive degrees, and a
certain endomorphism bialgebra, $\eA$, which coacts on all
components of $\K(\sA)$. These components therefore define elements
of the representation ring $R_{\eA}$ of $\eA$, and exactness of
$\K(\sA)$ in positive degrees yields an equation in the power series
ring $R_{\eA}\llbracket t \rrbracket$. Due to the specific form of
$\K(\sA)$, which is constructed from $\sA$ together with its
so-called dual algebra $\sA^!$, the equation in question does indeed
state that $\rho_1\cdot \rho_2 = 1$ holds for suitable $\rho_1,
\rho_2 \in R_{\eA}\llbracket t \rrbracket$. The last step in
deriving a MT for $\sA$ consists in using (super-)characters to
transport the abstract duality equation $\rho_1\cdot \rho_2 = 1$
from $R_{\eA}\llbracket t \rrbracket$ to the power series ring over
the algebra ${\eA}$, where it takes a more explicit and useable
form. Here then is the flow chart of our approach:
\begin{center}
$$
\xymatrix{
    *\txt{$N$-Koszul algebra \phantom{} \\ $\sA$}\ar@{->}[r]
    & *\txt{\ exact Koszul complex \phantom{} \\ $\K(\sA)$ } \ar@{->}[r]
    & *\txt{\ duality equation\phantom{} \\ in $R_{\eA}\llbracket t \rrbracket$}
    \ar@{->}[r] & *\txt{\phantom{x} MT for $\sA$}
    }
$$
\end{center}
The actual labor involved in this process consists in the explicit
evaluation of (super-)characters at the last arrow above. This step
is often facilitated by specializing the bialgebra $\eA$, which is
highly noncommutative, to a more familiar algebra $\sB$ via a
homomorphism $\eA \to \sB$. For example:
\begin{itemize}
\item
MacMahon's original MT \cite{paMM60} follows in the manner described
above by starting with $\sA = \sO(\k^d) = \k[x_1,\dots,x_d]$, the
ordinary polynomial algebra or ``affine space", and restricting the
resulting MT over $\en \sO(\k^d)$ to the coordinate ring of $d
\times d$-matrices, $\sO(\Mat_d(\k)) = \k[x^i_j \mid 1 \le i,j \le
d]$.
\item
As was explained in \cite{HL}, taking ``quantum affine space"
$\sO_\q(\k^d)$ as the point of departure one arrives at the quantum
MT of Garoufalidis, L{\^e} and Zeilberger \cite{GLZxx} (and
Konvalinka and Pak \cite{mKiP} in the multi-parameter case). The
endomorphism bialgebra of $\sO_\q(\k^d)$ is exactly the algebra of
right-quantum matrices as defined in \cite{GLZxx}.
\item
Berger's $N$-symmetric algebra \cite{rB01} leads to the
$N$-generalization of the MT proved by Etingof and Pak \cite{EP06}
using the above approach, again after restricting to
$\sO(\Mat_d(\k))$.
\end{itemize}

\subsection{} \label{SS:this article}

The present article aims to set forth an extension of the existing
theory of $N$-homogeneous algebras to the category $\Vect_{\k}^s$ of
vector superspaces over some base field $\k$. While this does not
give rise to principal obstacles given that \cite{yM88} and
\cite{rBmDVmW} are at hand as guiding references, the setting of
superalgebra requires careful consideration of the order of terms
and the so-called ``rule of signs" will be ubiquitous in our
formul{\ae}. In view of the potential interdisciplinary interest of
this material, we have opted to keep our presentation
reasonably self-contained and complete. 

Therefore, in Sections~\ref{S:linear} and \ref{S:supercharacter}, we
deploy the requisite background material from superalgebra in some
detail before turning to $N$-homogeneous superalgebras in
Section~\ref{S:homog}. The latter section, while following the
general outline of \cite{yM88} and \cite{rBmDVmW} rather closely,
also offers explicit discussions of a number of important examples.
We interpolate the pure even and pure odd Yang-Mills algebras
defined by Connes and Dubois-Violette \cite{aCmDV02}, \cite{aCmDV03}
by a family of superalgebras $\YM^{p|q}$ and give a unified
treatment of these algebras. (It turns out, however, that the mixed
algebras $\YM^{p|q}$, with $p$ and $q$ both nonzero, are less
well-behaved than the pure cases.) Moreover, we discuss a superized
version of the $N$-symmetric algebras of Berger \cite{rB01}.
Finally, in Example~\ref{EX:SRN}, we introduce new $N$-homogeneous
superversions of the symmetric algebra and the Grassmann algebra of
a vector superspace $V$; these are associated with any Hecke
operator $\hR \colon V^{\otimes 2} \to V^{\otimes 2}$ and will be
denoted by $\Sym_{\hR,N}$ and $\Gr_{\hR,N}$, respectively.

Sections~\ref{S:koszul} and \ref{S:duality} contain our main
results: Theorem~\ref{T:HeckeKoszul} shows that the superalgebras
$\Sym_{\hR,N}$ and $\Gr_{\hR,N}$ are in fact $N$-Koszul, and
Theorem~\ref{T:MTSN} is superized version of the aforementioned
$N$-generalized MT of Etingof and Pak \cite[Theorem 2]{EP06}. The
special case $N=2$ of Theorem~\ref{T:MTSN} is a superization of the
original master theorem of MacMahon \cite{paMM60}. The present
article was motivated in part by a comment in Konvalinka and Pak
\cite[13.4]{mKiP} asking for a ``real" super-analog of the classical
MT.

\subsection{} \label{SS:recent}

A considerable amount of research has been done by mathematical
physicists on various quantum matrix identities.
Some of these investigations have been carried out in a super
setting; see, e.g., Gurevich,  Pyatov and Saponov \cite{GPS05},
\cite{GPS06} and the references therein. However, the techniques
employed in these articles appear to be quite different from ours.

After submitting this article, we also learned of recent work of
Konvalinka \cite{mK07a}, \cite{mK07b} which not only concerns
MacMahon's MT but also other matrix identities such as the
determinantal identity of Sylvester. These identities are proved in
\cite{mK07a}, \cite{mK07b} by combinatorial means in various
noncommutative settings including the right-quantum matrix algebra
$\en{\sO_q(\k^d)}$.

\subsection{} \label{SS:generalized}

We conclude this Introduction by reviewing the precise definitions
of $N$-homogeneous and $N$-Koszul algebras. Our basic reference is
Berger \cite{rB01}; see also \cite{BGS96}, \cite{rBmDVmW},
\cite{eGetal04}.

Let $\sA$ be a connected $\Z_{\ge 0}$-graded algebra over a field
$\k$; so $\sA = \bigoplus_{n \ge 0} \sA_n$ for $\k$-subspaces
$\sA_n$ with $\sA_0 = \k$ and $\sA_n \sA_{m} \subseteq \sA_{n+m}$.
Choose a minimal generating set for the algebra $\sA$ consisting of
homogeneous elements of positive degree; this amounts to choosing a
graded basis for a graded subspace $V \subseteq \sA_+ = \bigoplus_{n
> 0} \sA_n$ such that $\sA_+ = \sA_+^2 \oplus V$. The grading of $V$
imparts a grading to the tensor algebra $\T(V)$ of the space $V$,
and we have a graded presentation
\begin{equation*}
\T(V)/I \iso \sA
\end{equation*}
for some graded ideal $I$ of $\T(V)$, the ideal of \emph{relations}
of $\sA$.

Recall that a graded vector space $M = \bigoplus_{n \in \Z} M_n$ is
said to live in degrees $\ge n_0$ if $M_n = 0$ for all $n < n_0$.
Note that the relation ideal $I$ lives in degrees $\ge 2$, because
$\T(V)_0 \oplus \T(V)_1 \subseteq \k \oplus V$ and $\k \oplus V$
injects into $\sA$. Fix an integer $N \ge 2$ and define the jump
function
\begin{equation} \label{E:jump1}
\ju_N(i) = \begin{cases} \tfrac{i}{2}N \quad &\text{if $i$ is even}\\
\tfrac{i-1}{2}N + 1 \quad &\text{if $i$ is odd}
\end{cases}
\end{equation}
The following proposition is identical with \cite[Proposition
2.1]{rBnM06} except for the fact that we do not a priori assume
$\sA$ to be generated in degree $1$. A proof is given in the
Appendix.

\begin{prop} \label{P:intro}
The ideal $I$ of relations of $\sA$ lives in degrees $\ge N$ if and
only if $\Tor^{\sA}_i(\k,\k)$ lives in degrees $\ge \ju_N(i)$ for
all $i \ge 0$.
\end{prop}

Following Berger \cite{rB01}, the graded algebra $\sA$ is said to be
\emph{$N$-Koszul} if $\Tor^{\sA}_i(\k,\k)$ is concentrated in degree
$\ju_N(i)$ for all $i \ge 0$. This implies that the space of algebra
generators $V$ is concentrated in degree $\ju_N(1) = 1$; so the
algebra $\sA$ is \emph{$1$-generated}. Moreover, choosing a minimal
set of homogeneous ideal generators for the relation ideal $I$
amounts to choosing a graded basis for a graded subspace $R
\subseteq I$ such that
\begin{equation} \label{E:I}
I = R \oplus \left(V \otimes I + I \otimes V\right)
\end{equation}
Then $\Tor^{\sA}_2(\k,\k) \cong R$ and so $R$ must be concentrated
in degree $\ju_N(2) = N$ when $\sA$ is $N$-Koszul. To summarize, all
$N$-Koszul algebras are necessarily $1$-generated and they have
defining relations in degree $N$; so there is a graded isomorphism
$$
\sA \cong \T(V)/(R) \qquad \text{with} \quad R \subseteq V^{\otimes
N}
$$
Such algebras are called \emph{$N$-homogeneous}.

We remark that Green et al.~\cite{eGetal04} have studied $N$-Koszul
algebras in the more general context where the grading $\sA =
\bigoplus_{n \ge 0} \sA_n$ is not necessarily connected ($\sA_0 =
\k$). In \cite[Theorem 4.1]{eGetal04}, it is shown that an
$N$-homogeneous algebra $\sA$ with $\sA_0$ split semisimple over
$\k$ is $N$-Koszul if and only if the Yoneda $\Ext$-algebra $E(\sA)
= \bigoplus_{n \ge 0} \Ext^n_{\sA}(\sA_0,\sA_0)$ is generated in
degrees $\le 2$.

Any $N$-homogeneous algebra $\sA$ whose generating space $V$ carries
a $\Z_2$-grading  and whose defining relations $R$ are $\Z_2$-graded
is naturally a $\k$-superalgebra, that is, $\sA$ has a
$\Z_2$-grading (``parity") besides the basic $\Z_{\ge 0}$-grading
(``degree"). As will be reviewed below, this extra structure
provides us with additional functions on Grothendieck rings, namely
superdimension and supercharacters, which lead to natural
formulations of the MT in a superized context. Note, however, that
the defining property of $N$-Koszul algebras makes no reference to
the $\Z_2$-grading of $\sA$. Thus, an $N$-homogeneous superalgebra
is Koszul precisely if it is Koszul as an ordinary $N$-homogeneous
algebra (forgetting the $\Z_2$-grading).

\subsection{} \label{SS:notations}

Throughout $\k$ is a commutative field and $\otimes$ stands for
$\otimes_{\k}$. Scalar multiplication in $\k$-vector spaces will
often, but not always, be written on the right while linear maps
will act from the left. We tacitly assume throughout that $\ch \k
\neq 2$; further restrictions on the characteristic of $\k$ will be
stated when required.


\section{Review of linear superalgebra} \label{S:linear}

\subsection{Vector superspaces} \label{SS:superspace}
A vector superspace over $\k$ is a $\k$-vector space $V$ equipped
with a grading by the group $\Z_2 = \Z/2\Z = \{\bar 0,\bar 1\}$.
Thus, we have a decomposition  $V = V_{\bar 0}\oplus V_{\bar 1}$
with $\k$-subspaces $V_{\bar 0}$ and $V_{\bar 1}$ whose elements are
called \emph{even} and \emph{odd}, respectively. In general, the
$\Z_2$-degree of a homogeneous element $a \in V$ is also called its
\emph{parity}; it will be denoted by $\dg{a} \in \Z_2$. Vector
superspaces over $\k$ form a category $\Vect_{\k}^s$ whose morphisms
are given by the linear maps preserving the $\Z_2$-grading; such
maps are also called \emph{even} linear maps.

The \emph{dimension} of an object $V$ of $\Vect_{\k}^s$ is the usual
$\k$-linear dimension. We shall use the notation
\begin{equation*}
d = \dim_\k V\,,\ p = \dim_\k V_{\bar 0}\quad \text{and} \quad q =
\dim_\k V_{\bar 1}
\end{equation*}
So $d = p+q$.  The \emph{superdimension} of a vector superspace $V$
with $d < \infty$ is defined by
\begin{equation*}\label{E:sdim}
\sdim V = p - q \in \Z
\end{equation*}
When working with a fixed basis $\{x_i\}$ of a given $V$ in
$\Vect_{\k}^s$ we shall assume that each $x_i$ is homogeneous; the
parity of $x_i$ will be denoted by $\dg{i}$. The basis
$x_1,x_2\dots$ is called \emph{standard} if $\dg{i}= \bar 0 \ (i \le
p)$ and $\dg{i} = \bar 1 \ (i>p)$.

\subsection{Tensors} \label{SS:tensors}
The tensor product $U\otimes V$ of vector superspaces $U$ and $V$ in
$\Vect_{\k}^s$ is the usual tensor product over $\k$ of the
underlying vector spaces equipped with the natural $\Z_2$-grading:
if $a,b$ are homogeneous elements then the parity of $a\otimes b$ is
$\dg{a} + \dg{b} \in \Z_2$. Instead of the usual symmetry
isomorphism $U\otimes V \iso V \otimes U$ for interchanging terms in
a tensor product we shall use the \emph{rule of signs}, that is, the
following functorial \emph{supersymmetry} isomorphism in
$\Vect_{\k}^s$:
\begin{equation}\label{E:symm}
c_{U,V} \colon U\otimes V \iso V\otimes U\ , \quad u\otimes v\mapsto
(-1)^{\dg{u}\dg{v}}v\otimes u 
\end{equation}
for $u,v$ homogeneous. (All formulas stated for homogeneous elements
only are to be extended to arbitrary elements by linearity.) The
supersymmetry isomorphisms $c_{U,V}$ satisfy $c_{V,U} \circ c_{U,V}
= \Id_{U \otimes V}$, and they are compatible with the usual
associativity isomorphims $a_{U,V,W} \colon (U \otimes V) \otimes W
\cong U \otimes (V \otimes W)$ in $\Vect_{\k}^s$, that is, they
satisfy the ``Hexagon Axiom"; see \cite[Def.~XIII.1.1]{cK95}.
Therefore, $\Vect_{\k}^s$ is a symmetric tensor category; the unit
object is the field $\k$, with parity $\bar 0$. See
\cite[Chap.~XIII]{cK95} or \cite{pDjM} for background on tensor
categories.

\subsection{Homomorphisms} \label{SS:hom}
The space $\Hom_\k(V,U)$ of all $\k$-linear maps between vector
superspaces $V$ and $U$ is again an object of $\Vect_{\k}^s$, with
grading $\Hom_\k(V,U)_{\bar 0} = \Hom_{\k}(V_{\bar 0},U_{\bar 0})
\oplus \Hom_\k(V_{\bar 1},U_{\bar 1})$ and $\Hom_{\k}(V,U)_{\bar 1}
= \Hom_\k(V_{\bar 0},U_{\bar 1}) \oplus \Hom_{\k}(V_{\bar 1},U_{\bar
0})$; so
$$
\Hom_\k(V,U)_{\bar 0} =
\Hom_{\Vect_{\k}^s}(V,U)
$$
In particular, the linear dual $V^* =
\Hom_{\k}(V,\k)$ belongs to $\Vect_{\k}^s$. Given homogeneous bases
$\{x_j\}$ of $V$ and $\{y_i\}$ of $U$ we can describe any $f \in
\Hom_\k(V,U)$ by its matrix $F = (F^i_j)$:
\begin{equation} \label{E:matrix}
f(x_j) = \sum_i y_i F^i_j
\end{equation}
When $f$ is an even map then  $F^i_j = 0$ unless $\dg{i} + \dg{j} =
\bar 0$.

For finite-dimensional vector superspaces, we have the following
functorial isomorphisms in $\Vect_{\k}^s$ (see, e.g.,
\cite[I.8]{gT04}):
\begin{equation} \label{E:iso1}
U \otimes V^* \cong \Hom_{\k}(V,U)
\end{equation}
via $(u\otimes f)(v) = u \langle f,v \rangle$, and
\begin{equation}  \label{E:iso2}
V_1^* \otimes \ldots \otimes  V_m^* \cong (V_m \otimes \ldots
\otimes V_1)^*
\end{equation}
via $\langle f_1 \otimes \ldots \otimes f_m,v_m \otimes \ldots
\otimes v_1\rangle = \prod_i \langle f_i,v_i\rangle$.
Here, we use the notation $\langle f,v \rangle = f(v)$ for the
evaluation pairing
\begin{equation*} \label{E:ev}
\ev_V = \langle\,.\,,.\,\rangle \colon V^* \otimes V \to \k
\end{equation*}
in $\Vect_{\k}^s$. Similarly, we have a pairing
\begin{equation*} \label{E:ev2}
V \otimes V^*
\stackrel{c_{V,V^*}}{\tto} V^* \otimes V \stackrel{\ev_V}{\tto} \k
\end{equation*}
which yields an isomorphism
\begin{equation}\label{E:**}
V \iso V^{**}
\end{equation}
in $\Vect_{\k}^s$.

The isomorphism \eqref{E:iso1} (which is valid as long as one of $U$
or $V$ is finite-dimensional) has the following explicit
description. Fix homogeneous bases $\{x_j\}$ of $V$ and $\{y_i\}$ of
$U$ and let $F = (F^i_j)$ be the matrix of a given $f \in
\Hom_\k(V,U)$ with respect to these bases, as in \eqref{E:matrix}.
Let $\{x^j\}$ be the dual basis of $V^*$, defined by $\langle
x^j,x_\ell\rangle = \delta^{j}_{\ell}$ (Kronecker delta). Then the
image of $f$ in $U \otimes V^*$ is given by $\sum_{i,j} y_i\otimes
x^j F^i_j$. Note also that $x_i$ and $x^i$ have the same parity.

Finally, if $U$, $V$ and $W$ are vector superspaces, with $U$
finite-dimensional, then the isomorphism $\Id \otimes c_{W,U^*}
\colon  V \otimes W \otimes U^* {\iso} V \otimes U^* \otimes W$
together with \eqref{E:iso1} yields an isomorphism
\begin{equation} \label{E:hom}
\Hom_{\k}(U,V \otimes W) \iso \Hom_{\k}(U,V) \otimes W
\end{equation}
in $\Vect_{\k}^s$ which is explicitly given by  $(f \otimes w)(u) =
(-1)^{\dg{w} \dg{u}} f(u) \otimes w$. Similarly, for vector
superspaces $U$, $U'$, $V$, $V'$ with $U$, $U'$ finite-dimensional,
there is an isomorphism
\begin{equation} \label{E:hom2}
\Hom_{\k}(U \otimes U',V \otimes V') \iso \Hom_{\k}(U,V) \otimes
\Hom_{\k}(U',V')
\end{equation}
in $\Vect_{\k}^s$ given by  $(f \otimes g)(u \otimes v) =
(-1)^{\dg{g} \dg{u}} f(u) \otimes g(v)$.

\subsection{Supertrace} \label{SS:str}
Let $V$ be a finite-dimensional object of $\Vect_{\k}^s$. The
\emph{supertrace} is the map
\begin{equation}\label{E:str}
\str \colon \End_\k(V) \underset{\eqref{E:iso1}}\iso V \otimes V^*
\underset{\eqref{E:ev2}}{\tto} \k
\end{equation}
in $\Vect_{\k}^s$. In order to describe the supertrace in terms of
matrices, fix a basis $\{ x_i \}$ of $V$ consisting of homogeneous
elements and let $F = (F^i_j)$ be the matrix of $f \in \End_{\k}(V)$
as in \eqref{E:matrix}. Then
\begin{equation*} \label{E:str2}
\str(f) = \sum_i (-1)^{\dg{i}} F^i_i
\end{equation*}
where $\dg{i}$ is the parity of $x_i$ (and of the dual basis vector
$x^i \in V^*$) as in \S\ref{SS:superspace}. Thus,
\begin{equation*}
\str(\Id_V) = \sdim V . 1_\k
\end{equation*}

\subsection{Action of the symmetric group} \label{SS:action}
Given vector superspaces $V_1,\dots,V_n$, we can consider the
morphism
$$
c_{i} \colon V_1 \otimes \dots \otimes V_i \otimes V_ {i+1} \otimes
\dots \otimes V_n \tto V_1 \otimes \dots \otimes V_{i+1} \otimes V_
{i} \otimes \dots \otimes V_n
$$
in $\Vect_{\k}^s$ which interchanges the factors $V_i$ and $V_{i+1}$
via $c_{V_i,V_{i+1}}$ and is the identity on all other factors. More
generally, for any $\si \in \Sy_n$, the symmetric group  consisting
of all permutations of $\{1,2,\dots,n\}$, one can define a morphism
$$
c_\sigma \colon V_1 \otimes \dots \otimes V_n \tto
V_{\sigma^{-1}(1)} \otimes \dots \otimes V_{\sigma^{-1}(n)}
$$
in $\Vect_{\k}^s$ as follows. Recall that $\Sy_n$ is generated by
the transpositions $\si_1,\dots,\si_{n-1}$ where $\si_i$
interchanges $i$ and $i+1$ and leaves all other elements of
$\{1,2,\dots,n\}$ fixed. The minimal length of a product in the
$\si_i$'s which expresses a given element $\si \in \Sy_n$ is called
the length of $\sigma$ and denoted $\ell(\si)$; it is given by
$$
\ell(\si) = \# \inv(\si) \quad\text{with}\quad \inv(\si) = \{ (i,j)
\mid i < j \text{ but } \si(i) > \si(j) \}
$$
Writing $\si \in \Sy_n$ as a product of certain $\si_i$, the
analogous product of the maps $c_{i}$ yields a morphism $c_\sigma$
as above. This morphism is independent of the way $\sigma$ is
expressed in terms of the transpositions $\si_i$; see
\cite[I.4.13]{gT04} or \cite[Theorem~XIII.1.3]{cK95}. If all $v_i
\in V_i$ are homogeneous then
\begin{equation} \label{E:c}
c_\si(v_1 \otimes \dots \otimes v_n) = (-1)^{\sum_{(i,j) \in
\inv(\si)} \dg{v_i}\dg{v_j}} v_{\si^{-1}(1)}  \otimes \dots \otimes
v_{\si^{-1}(n)}
\end{equation}
For example, if all $v_i$ are even then the $\pm$-sign on the right
is $+$, and if all $v_i$ are odd then it is $\sign(\si)$, the
signature of $\si$.

Taking all $V_i = V$ we obtain a representation $c \colon \Sy_n \tto
\Aut_{\Vect_{\k}^s}(V^{\otimes n})$ where $V^{\otimes n} = V \otimes
\dots \otimes V$ ($n$ factors). Letting $\k[\Sy_n]$ denote the group
algebra of the symmetric group, this extends uniquely to an algebra
map
\begin{equation} \label{E:rep}
c \colon \k[\Sy_n] \tto \End_{\Vect_{\k}^s}(V^{\otimes n})
\end{equation}
We will write $c_a:= c(a)$ for $a \in \k[\Sy_N]$.

For the dual superspace $V^*$, besides the above representation $c
\colon \k[\Sy_n] \tto \End_{\Vect_{\k}^s}(V^{*\otimes n})$, we also
have the \emph{contragredient representation}
\begin{equation*} \label{E:contra}
c^* \colon \k[\Sy_n] \tto \End_{\Vect_{\k}^s}(V^{*\otimes n})
\end{equation*}
for the pairing $\langle\,.\,,.\,\rangle \colon V^{*\otimes n}
\otimes V^{\otimes n} \to \k$ in \eqref{E:iso2}. Explicitly,
$$
\langle c^*_a(x),y \rangle = \langle x, c_{a^*}(y) \rangle
$$
for all $a \in \k[\Sy_n]$, $x \in V^{*\otimes n}$ and $y \in
V^{\otimes n}$. Here, $\,.\,^* \colon \k[\Sy_n] \to \k[\Sy_n]$ is
the involution sending $\sigma \in \Sy_n$ to $\sigma^{-1}$. These
two representations are related by
\begin{equation} \label{E:contra2}
c^*_a = c_{\tau a \tau}
\end{equation}
where $\tau = (1,n)(2,n-1)\dots \in \Sy_n$ is the order reversal
involution. One only needs to check \eqref{E:contra2} for the
transpositions $a = \si_i$, which is straightforward.

\subsection{Hecke algebras} \label{SS:Hecke}
We recall some standard facts concerning Hecke algebras; these are
suitable deformations of the group algebra $\k[\Sy_n]$ considered
above. For details, see \cite{rDgJ86}, \cite{rDgJ91}.

Fix $0 \neq q \in \k$. The Hecke algebra $\He_{n,q}$ is generated as
$\k$-algebra by elements $T_1,\dots, T_{n-1}$ subject to the
relations
\begin{equation} \label{E:Hecke1}
\begin{aligned}
&(T_i + 1)(T_i - q) = 0  \\
&T_iT_{i+1}T_i = T_{i+1}T_iT_{i+1} \\
&T_iT_j = T_jT_i \quad \text{if $|i-j| \ge 2$}
\end{aligned}
\end{equation}
When $q=1$, one has an isomorphism $\He_{n,1} \iso
\k[\Sy_n]$, $T_i \mapsto \si_i$ where $\si_i$ is the transposition
$(i,i+1)$ as in \S\ref{SS:action}. The algebra $\He_{n,q}$ has a
$\k$-basis $\{ T_\si \mid \si \in \Sy_n \}$ so that
\begin{enumerate}
\item[(i)] $T_{\Id} = 1$ and $T_{\si_i} = T_i$;
\item[(ii)] $T_\si T_{\si_i} =
\begin{cases} T_{\si\si_i} &\text{if $\ell(\si\si_i) = \ell(\si) + 1$}; \\
q T_{\si\si_i} + (q-1) T_\si &\text{otherwise}
\end{cases}
$
\end{enumerate}
By $\k$-linear extension of the rule
\begin{equation*} \label{E:Heckeinvol}
T_\si^*:= T_{\si^{-1}}\qquad (\si \in \Sy_n)
\end{equation*}
one obtains an involution $\,.\,^* \colon \He_{n,q} \to \He_{n,q}$.
Moreover, the elements $T_i':= -q T_i^{-1} = q - 1 - T_i$ also
satisfy relations \eqref{E:Hecke1}. Therefore,
\begin{equation} \label{E:Heckeaut}
\alpha(T_i ):= -q T_i^{-1}
\end{equation}
defines an algebra automorphism $\alpha \colon \He_{n,q} \to
\He_{n,q}$ of order $2$.

The Hecke algebra $\He_{n,q}$ is always a symmetric algebra, and
$\He_{n,q}$ is a split semisimple $\k$-algebra iff the following
condition is satisfied:
\begin{equation} \label{E:Heckess}
[n]_q!:= \prod_{i=1}^n [i]_q \neq 0 \quad \text{ where $[i]_q :=
1+q+\dots+q^{i-1}$}
\end{equation}
More precisely, if \eqref{E:Heckess} holds then
\begin{equation} \label{E:Heckess2}
\He_{n,q} \cong \bigoplus_{\lambda\vdash n} \Mat_{d_\lambda \times
d_\lambda}(\k)
\end{equation}
where $\lambda$ runs over all partitions of $n$ and $d_\lambda$
denotes the number of standard $\lambda$-tableaux. The only
partitions $\lambda$ with $d_\lambda = 1$ are $\lambda = (n)$ and
$\lambda = (1^n)$. The central primitive idempotents of $\He_{n,q}$
for these partitions are given by
\begin{equation} \label{E:Heckeid1}
X_n := \frac{1}{[n]_q!}\sum_{\si \in \Sy_n} T_\si
\end{equation}
and
\begin{equation} \label{E:Heckeid2}
Y_n := \frac{1}{[n]_{q^{-1}}!}\sum_{\si \in \Sy_n} (-q)^{-\ell(\si)}
T_\si
\end{equation}
These idempotents are usually called the \emph{$q$-symmetrizer} and
the \emph{$q$-antisymme\-trizer}, respectively. One has
\begin{equation} \label{E:Heckeid3}
X_n T_\si = T_\si X_n = q^{\ell(\si)}X_n \quad \text{and} \quad Y_n
T_\si = T_\si Y_n = (-1)^{\ell(\si)}Y_n
\end{equation}
for $\si \in \Sy_n$. Furthermore, $\alpha(X_n) = Y_n$.

For later use, we note the following well-known consequence of
\eqref{E:Heckeid3}. If $M$ is any $\He_{n,q}$-module, with
corresponding representation $\mu \colon  \He_{n,q} \to \End_\k(M)$,
then
\begin{equation} \label{E:imXn}
\im(\mu(X_n)) = \bigcap_{i=1}^{n-1} \im(\mu(T_i)+1)
\end{equation}
Indeed, \eqref{E:Heckeid3} implies that $X_n = [2]_q^{-1}(T_i +
1)X_n$, which yields the inclusion $\subseteq$. On the other hand,
any $m \in \bigcap_{i=1}^{n-1} \im(\mu(T_i)+1)$ satisfies $(\mu(T_i)
- q)(m) = 0$ for all $i$, by \eqref{E:Hecke1}. Therefore,
$\mu(T_\si)(m) = q^{\ell(\si)}m$ holds for all $\si \in \Sy_n$, and
hence $\mu(X_n)(m) = \frac{1}{[n]_q!}\sum_{\si \in \Sy_n}
q^{\ell(\si)}m = m$. This proves $\supseteq$.

\subsection{Hecke operators} \label{SS:Heckeop}
Again, let $0 \neq q \in \k$. A Hecke operator (associated to $q$)
on a vector superspace $V$ is a morphism $\hR \colon V^{\otimes 2}
\to V^{\otimes 2}$ in $\Vect_{\k}^s$ satisfying the Hecke equation
\begin{equation*}  \label{E:hecke}
(\hR+1)(\hR-q) = 0
\end{equation*}
and the Yang-Baxter equation 
\begin{equation*}  \label{E:braid}
 \hR_1\hR_2\hR_1 = \hR_2\hR_1\hR_2
\end{equation*}
where $\hR_1 := \hR \otimes \Id_V \colon V^{\otimes 3} \to
V^{\otimes 3}$ and similarly $\hR_2 := \Id_V \otimes \hR$.

The Hecke equation implies that $\hR$ is invertible. Moreover, if
$\hR$ is a Hecke operator associated to $q$ then so is $-q\hR^{-1}$.

Defining $\rho(T_i) := \Id_V^{\otimes i-1}\otimes \hR \otimes
\Id_V^{\otimes n- i-1}$, one obtains a representation
\begin{equation} \label{E:Heckerep}
\rho = \rho_{n,\hR} \colon \He_{n,q} \tto
\End_{\Vect_{\k}^s}(V^{\otimes n})
\end{equation}
The representations $\rho_{n,\hR}$ and $\rho_{n,-q\hR^{-1}}$ are
related by $\rho_{n,-q\hR^{-1}} = \rho_{n,\hR} \circ \alpha$, where
$\alpha$ is the automorphism of $\He_{n,q}$ defined in
\eqref{E:Heckeaut}.

\begin{example} \label{EX:Heckec}
The supersymmetry operator $c_{V,V} \colon V^{\otimes 2} \to
V^{\otimes 2}$ in \eqref{E:symm} is a Hecke operator associated to
$q=1$, as is its negative, $-c_{V,V}$. The representation
$\rho_{c_{V,V}}$ of $\He_{n,1} = \k[\Sy_n]$ in \eqref{E:Heckerep} is
identical with \eqref{E:rep}.
\end{example}

\begin{example}[superized Drinfel'd-Jimbo \cite{yM89}, \cite{phH02a}]
\label{EX:HeckDJ} Let $x_1,\dots,x_d$ be a standard basis of $V$ as
in \S\ref{SS:superspace}. The super analog $\hR = \hR^{DJ}$ of the
standard Drinfel'd-Jimbo Hecke operator is defined as follows.
Writing
\begin{equation*} \label{E:DJ1}
\hR(x_i \otimes x_j) = \sum_{k,l} x_k \otimes x_l \hR^{k,l}_{i,j}
\end{equation*}
the matrix components $\hR^{k,l}_{i,j} \in \k$ are given by
\begin{equation*} \label{E:DJ2}
\hR^{k,l}_{i,j} = \frac{q^2 - q^{2\e_{i,j}}}{1 + q^{2\e_{i,j}}}
\delta^{k,l}_{i,j} + (-1)^{\dg{i}\dg{j}} \frac{q^{\e_{i,j}}(q^2 +
1)}{1 + q^{2\e_{i,j}}} \delta^{l,k}_{i,j}
\end{equation*}
Here, $\e_{i,j} = \sign(i-j)$. Thus,
\begin{equation} \label{E:HeckDJ}
\begin{aligned}
\hR^{ii}_{ii} &= q^2 & &\text{if $\dg{i} = \bar 0$} \\
\hR^{ii}_{ii} &= -1 & &\text{if $\dg{i} = \bar 1$} \\
\hR^{ij}_{ij} &= q^2-1 & &\text{if $i < j$} \\
\hR^{ji}_{ij} &= (-1)^{\dg{i}\dg{j}} q & &\text{if $i \neq j$}
\end{aligned}
\end{equation}
and $\hR^{k,l}_{i,j} = 0$ in all other cases. One checks that $\hR$
is a Hecke operator that is associated to $q^2$.
\end{example}


\section{The supercharacter} \label{S:supercharacter}

\subsection{Superalgebras, supercoalgebras etc.} \label{SS:superalg}
An algebra $\sA$ in $\Vect_{\k}^s$ is called a \emph{superalgebra}
over $\k$; this is just an ordinary $\k$-algebra such that the unit
map $\k \to \sA$ and the multiplication
$$
\mu \colon \sA \otimes \sA
\to \sA
$$
are morphisms in $\Vect_{\k}^s$. In other words, $\sA$ is a
$\Z_2$-graded $\k$-algebra in the usual sense: $\sA = \sA_{\bar
0}\oplus \sA_{\bar 1}$ with $\k$-subspaces $\sA_{\bar 0}$ and
$\sA_{\bar 1}$ such that $\sA_{\bar r}\sA_{\bar s} \subseteq
\sA_{\overline{r+s}}$. Homomorphisms of superalgebras, by
definition, are algebra maps in $\Vect_{\k}^s$, that is, they
preserve the $\Z_2$-grading.

If $V$ is a vector superspace in $\Vect_{\k}^s$ then the tensor
algebra $\T(V) = \bigoplus_{n \ge 0} V^{\otimes n}$ is a
superalgebra via the $\Z_2$-grading of each $V^{\otimes n}$ as in
\S\ref{SS:tensors}. In general, if $\sA$ is any superalgebra, then
by selecting a $\Z_2$-graded subspace $V \subseteq \sA$ which
generates the algebra $\sA$, we obtain a canonical isomorphism of
superalgebras
\begin{equation} \label{E:superisom}
\T(V)/(R) \iso \sA
\end{equation}
where $(R)$ is the two-sided ideal of $\T(V)$ that is generated by a
$\Z_2$-graded linear subspace $R \subseteq \T(V)$.

Given superalgebras $\sA$ and $\sB$, the tensor product $\sA \otimes
\sB$ is the superalgebra with the usual additive structure and
grading and with multiplication $\mu_{\sA\otimes\sB}$ defined by
using the supersymmetry map \eqref{E:symm}: $\mu_{\sA\otimes\sB} =
(\mu_{\sA} \otimes \mu_{\sB})\circ(\Id_A \otimes c_{\sB,\sA} \otimes
\Id_{\sB})$ or, explicitly,
\begin{equation*} \label{E:tensmult}
(a \otimes b)(a' \otimes b') = (-1)^{\dg{a'}\dg{b}} aa' \otimes bb'
\end{equation*}
for homogeneous $a' \in \sA$ and $b \in \sB$. In other words, the
canonical images of $\sA$ and $\sB$ in $\sA\otimes\sB$
\emph{supercommute}, in the sense that the supercommutator
\begin{equation} \label{E:supercomm}
[a,b] = ab - (-1)^{\dg{a}\dg{b}} ba
\end{equation}
vanishes for any pair of homogeneous elements $a \in \sA$ and $b \in
\sB$.

Supercoalgebras, superbialgebras etc. are defined similarly as
suitable objects of $\Vect_{\k}^s$ such that all structure maps are
maps in $\Vect_{\k}^s$. The compatibility between the
comultiplication $\Delta$ and the multiplication of a superbialgebra
$\sB$ amounts to the following rule:
\begin{equation*} \label{E:supercomult}
\Delta(ab) =
\sum_{(a),(b)}(-1)^{\dg{a}_{(2)}\dg{b}_{(1)}}a_{(1)}b_{(1)} \otimes
a_{(2)}b_{(2)}
\end{equation*}
for homogeneous elements $a,b \in \sB$. Here we use the Sweedler
notation $\Delta(a) = \sum_{(a)} a_{(1)} \otimes a_{(2)}$ and
$a_{(1)}, a_{(2)}$ are chosen homogeneous with $\dg{a}_{(1)} +
\dg{a}_{(2)} = \dg{a}$.

\begin{example}[Symmetric superalgebra {\cite[3.2.5]{yM97}}] \label{EX:SV}
The symmetric superalgebra of a given $V$ in $\Vect_{\k}^s$ is
defined by
\begin{equation*}
\Sym(V) = \T(V)/\left( [v,w]_\otimes \mid v,w \in V  \right)
\end{equation*}
where $[v,w]_\otimes$ is the supercommutator \eqref{E:supercomm} in
$\T(V)$. Ignoring parity, $\Sym(V)$ is isomorphic to $\Sym(V_{\bar
0}) \otimes \Gr(V_{\bar 1})$, where $\Sym(\,.\,)$ and $\Gr(\,.\,)$
denote the ordinary symmetric and exterior (Grassmann) algebras,
respectively. The symmetric superalgebra is a Hopf superalgebra:
comultiplication $\Delta \colon \Sym(V) \to \Sym(V) \otimes \Sym(V)$
is given by $\Delta(v) = v \otimes 1 + 1 \otimes v$ for $v \in V$
and extension to all of $\Sym(V)$ by multiplicativity. Similarly,
the counit $\e \colon \Sym(V) \to \k$ is given by $\e(v) = 0$ and
the antipode $\sS \colon \Sym(V) \to \Sym(V)$ by $\sS(v) = -v$ for
$v \in V$.
\end{example}

\subsection{Comodules} \label{SS:comod}

We refer to \cite[Chap.~III]{cK95} for background on comodules,
comodule algebras etc.

Given a superbialgebra $\sB$, we let $\comod_{\sB}$ denote the
category of all right $\sB$-comodules and $\sB$-comodule maps in
$\Vect_{\k}^s$. Thus, for any object $V$ in $\comod_{\sB}$, we have
a ``coaction" morphism
$$
\delta_V \colon V \to V \otimes \sB
$$
in $\Vect_{\k}^s$. If $x_1,\dots,x_d$ is a fixed basis of $V$
consisting of homogeneous elements, with $\dg{i}$ denoting the
parity of $x_i$ as before, then we will write
\begin{equation} \label{E:deltamatrix}
\delta_V(x_j) = \sum_i x_i \otimes b^i_j \quad\text{with}\quad b^i_j
\in B_{\dg{i} + \dg{j}}
\end{equation}

The tensor product of vector superspaces makes $\comod_{\sB}$ into a
tensor category: if $U$ and $V$ are in $\comod_{\sB}$ then $\sB$
coacts on $U \otimes V$ by
\begin{equation} \label{E:deltatensor}
\delta_{U\otimes V} \colon U\otimes V\stackrel{\delta_U \otimes
\delta_V}{\tto} U\otimes \sB\otimes V\otimes \sB
\stackrel{c_{B,V}}{\tto} U \otimes V \otimes \sB \otimes \sB
\stackrel{\Id \otimes \mu_{\sB}}\tto U\otimes V\otimes \sB
\end{equation}
If $\sB$ is supercommutative as a superalgebra then the
supersymmetry $c_{U,V}$ is a $\sB$-comodule morphism, i.e.,
$\delta_{V\otimes U} \circ c_{U,V} = \left(c_{U,V} \otimes
\Id_{\sB}\right) \circ \delta_{U\otimes V}$. Therefore
$\comod_{\sB}$ is a \emph{symmetric} tensor category in this case .

\subsection{The supercharacter map} \label{SS:sch}
Let $\sB$ denote a superbialgebra and let $V$ be a finite
dimensional object in $\comod_{\sB}$. The coaction $\delta_V$ is an
even map in $\Hom_{\k}(V,V \otimes \sB)$. Consider the following
morphism in $\Vect_{\k}^s$:
\begin{equation} \label{E:sch1}
\schi \colon \End_{\k}(V) \stackrel{\delta_V\circ(\,.\,)}{\tto}
\Hom_{\k}(V,V \otimes \sB) \underset{\eqref{E:hom}}{\iso}
\End_{\k}(V) \otimes \sB \stackrel{\str \otimes \Id}{\tto} \k
\otimes \sB = \sB
\end{equation}
where $\str$ is the supertrace as in \eqref{E:str}. This map will be
called the \emph{supercharacter} map of $V$. Forgetting parity and
viewing all elements as even, the supertrace becomes the ordinary
trace and the supercharacter becomes the usual character. These will
be denoted by $\tr$ and $\chi$, respectively.

In particular, we have the element
\begin{equation*} \label{E:sch2}
\schi_V := \schi(\Id_V) \in \sB_{\bar 0}
\end{equation*}
To obtain explicit formulas, fix a basis $x_1,\dots,x_d$ of $V$
consisting of homogeneous elements and let $(F^i_j)$ and $(b^i_j)$
be the matrices of $f \in \End_{\k}(V)$ and of $\delta_V$ with
respect to this basis as in \eqref{E:matrix}, \eqref{E:deltamatrix}.
Then
\begin{equation} \label{E:sch3}
\schi(f) = \sum_{i,j} (-1)^{\dg{i} \dg{j}} b^i_jF^j_i
\end{equation}

Let $\e \colon \sB \to \k$ denote the counit of $\sB$. Then $x_j =
\sum_i x_i  \e(b^i_j)$ holds in \eqref{E:deltamatrix}. Hence
$\e(b^i_j) = \delta^i_j.1_\k$ and \eqref{E:sch3} gives
\begin{equation} \label{E:sch3a}
\e (\schi(f)) = \str(f)
\end{equation}

When $f$ is even formula \eqref{E:sch3} becomes $\schi(f) =
\sum_{i,j} (-1)^{\dg{i}} b^i_jF^j_i$, because $F^j_i = 0$ unless
$\dg{i} + \dg{j} = \bar 0$. In particular,
\begin{equation} \label{E:sch4}
\schi_V = \sum_i (-1)^{\dg{i}} b^i_i
\end{equation}
In the following, we let $\com_{\sB}$ denote the full subcategory of
$\comod_{\sB}$ consisting of all objects that are finite-dimensional
over $\k$. The supercharacter has the following properties analogous
to standard properties of the ordinary character.

\begin{lem} \label{L:sch}
Let $\sB$ denote a superbialgebra and let $U$, $V$ and $W$ be
objects of $\com_{\sB}$.
\begin{enumerate}
\item
If $f \colon V \to U$ and $g \colon U \to V$ are $\sB$-comodule maps
(not necessarily even) then
$$\schi(f \circ g) = (-1)^{\dg{f}\dg{g}} \schi(g \circ f)$$

\item
For $f \in \End_{\k}(V)$, $g \in \End_{\k}(U)$ view $f\otimes g
\in \End_{\k}(V \otimes U)$ as in \eqref{E:hom2}. Then
$$\schi(f\otimes g) = \schi(f)\schi(g)$$

\item
Given an exact sequence $0 \to U \stackrel{\mu}{\tto} V
\stackrel{\mu}{\tto} W \to 0$ in $\com_{\sB}$, let $f \in
\End_{\k}(V)$ be such that $f(\mu(U)) \subseteq \mu(U)$, and let $g
\in \End_{\k}(U)$, $h \in \End_{\k}(W)$ be the maps induced by $f$.
Then
$$
\schi(f) = \schi(g) + \schi(h)
$$
In particular, $\schi_V = \schi_U + \schi_W$. Moreover, if $f \in
\End_{\com_{\sB}}(V)$ is a projection (i.e., $f^2 = f$) then
$\schi(f)=\schi_{\im f}$.
\end{enumerate}
\end{lem}

\begin{proof}
(a) Let $T_V$ denote the map $\Hom_{\k}(V,V \otimes \sB) \tto \sB$
in \eqref{E:sch1}; so $\schi(f) = T_V(\delta_V \circ f)$. Since $f$
and $g$ are comodule maps, we have $\delta_{U} \circ f = (f \otimes
\Id_{\sB}) \circ \delta_V$ and similarly for $g$. Putting $h =
\delta_{U} \circ f \in \Hom_{\k}(V,U \otimes \sB)$ we obtain
$\schi(f \circ g) = T_{U}(\delta_{U} \circ f \circ g) = T_{U}(h
\circ g)$ and $\schi(g \circ f) = T_{V}(\delta_{V} \circ g \circ f)
= T_V((g \otimes \Id_{\sB})\circ h)$. Therefore, we must show that
$$
T_{U}(h \circ g) = (-1)^{\dg{f}\dg{g}} T_V((g \otimes
\Id_{\sB})\circ h)
$$
Using the identification $\Hom_{\k}(V,U \otimes \sB)  \cong
\Hom_{\k}(V,U)\otimes \sB$ as in \eqref{E:hom}, write $h = \sum_i
f_i \otimes b_i$ with $f_i \in \Hom_{\k}(V,U)$, $b_i \in \sB$, and
$\dg{f_i} + \dg{b_i} = \dg{h} = \dg{f}$. Then $h \circ g \in
\Hom_{\k}(U,U \otimes \sB)$ becomes the element $(\sum_i f_i \otimes
b_i) \circ g = \sum_i (-1)^{\dg{b_i}\dg{g}} (f_i \circ g) \otimes
b_i \in \End_{\k}(U) \otimes \sB$, and $(g \otimes \Id_{\sB})\circ h
\Hom_{\k}(V,V \otimes \sB)$ becomes $\sum_i (g \circ f_i) \otimes
b_i$. The standard identity $\str(f_i \circ g) =
(-1)^{\dg{f_i}\dg{g}} \str(g \circ f_i)$ (cf., e.g., \cite[p.~165 \S
3(b)]{yM97}) now yields
\begin{align*}
T_{U}(h \circ g) &= \sum_i (-1)^{\dg{b_i}\dg{g}} \str(f_i \circ g)
\otimes b_i  \\
&= \sum_i (-1)^{\dg{b_i}\dg{g}+ \dg{f_i}\dg{g}} \str(g \circ f_i)
\otimes b_i \\
&= (-1)^{\dg{f}\dg{g}} T_V((g \otimes \Id_{\sB})\circ h)
\end{align*}
as desired.

(b) Fix homogeneous $\k$-bases $\{x_i\}$ and $\{y_\ell\}$ of $V$ and
$U$, respectively, and write $\dg{x_i} = \dg{i}$, $\dg{y_\ell} =
\dg{\ell}$ as usual. Moreover, let $(F^i_j)$ and $(G^\ell_m)$ be the
matrices of $f$ and $g$ for these bases, as in \eqref{E:matrix}.
Then $\{x_i \otimes y_\ell\}$ is a basis of $V \otimes U$, with $x_i
\otimes y_\ell$ having parity $\dg{i} + \dg{\ell}$. Moreover,
\begin{align*}
(f \otimes g)(x_j \otimes y_m) &= (-1)^{\dg{g}\dg{j}} f(x_j) \otimes
g(y_m) \\
&= (-1)^{\dg{g}\dg{j}} \sum_i x_iF^i_j \otimes \sum_\ell
y_\ell G^\ell_m \\
&= \sum_{i,\ell} x_i \otimes y_\ell \Phi^{i,\ell}_{j,m} \qquad
\text{with $\Phi^{i,\ell}_{j,m} = (-1)^{(\dg{\ell} + \dg{m})\dg{j}}
F^i_j G^\ell_m$}
\end{align*}
because $G^m_\ell = 0$ unless $\dg{\ell} + \dg{m} = \dg{g}$.
Similarly, writing $\delta_V(x_j) = \sum_i x_i \otimes b^i_j$ with
$b^i_j \in \sB_{\dg{i} + \dg{j}}$ and $\delta_U(y_m) = \sum_\ell
y_\ell \otimes c^\ell_m$ with $c^\ell_m \in \sB_{\dg{\ell} +
\dg{m}}$, one obtains using \eqref{E:deltatensor}
$$
\delta_{V\otimes U}(x_j \otimes y_m) = \sum_{i,\ell} x_i \otimes
y_\ell \otimes \Psi^{i,\ell}_{j,m} \qquad \text{with
$\Psi^{i,\ell}_{j,m} = (-1)^{(\dg{i} + \dg{j})\dg{\ell}} b^i_j
c^\ell_m$}
$$
Therefore, formula \eqref{E:sch3} becomes
\begin{align*}
\schi(f \otimes g) &= \sum_{i,\ell,j,m} (-1)^{(\dg{i} +
\dg{\ell})(\dg{j} + \dg{m})} \Psi^{i,\ell}_{j,m} \Phi^{j,m}_{i,\ell}
\\
&= \sum_{i,\ell,j,m} (-1)^{\dg{i}\dg{j} + \dg{\ell}\dg{m}} b^i_j
F^j_i c^\ell_m G^m_\ell \\
&= \schi(f) \schi(g)
\end{align*}

(c) Choose a basis $\{ x_i \}$ of $V$ consisting of homogeneous
elements so that $x_i = \mu(y_i)$ for $i \le \dim U$ and let
$(F^i_j)$ be the matrix of $f$ for this basis. Then $F^i_j = 0$ for
$i > \dim U$, $j \le \dim U$. Moreover, the $y_i$ form a basis of
$U$ and the $z_i = \pi(x_i)$ form a basis of $W$, and the matrices
of $g$ and $h$ for these bases are $(F^i_j)_{i,j \le \dim U}$ and
$(F^i_j)_{i,j
> \dim U}$, respectively. Similarly, if $(b^i_j)$ is the matrix of $\delta_V$
with respect to the basis basis $\{ x_i \}$ as in
\eqref{E:deltamatrix} then $b^i_j = 0$ for $i > \dim U$, $j \le \dim
U$, and the matrices of $\delta_U$ and $\delta_W$ for the given
bases are $(b^i_j)_{i,j \le \dim U}$ and $(b^i_j)_{i,j
> \dim U}$, respectively. Therefore,
\begin{align*}
\schi(f) &= \sum_{i,j} (-1)^{\dg{i} \dg{j}} b^i_jF^j_i \\
&= \sum_{i,j \le \dim U} (-1)^{\dg{i} \dg{j}} b^i_jF^j_i + \sum_{i,j
> \dim U} (-1)^{\dg{i} \dg{j}} b^i_jF^j_i \\
&= \schi(g) + \schi(h)
\end{align*}
The remaining assertions are clear.
\end{proof}

\subsection{The Grothendieck ring} \label{SS:grothendieck}

Let $\sB$ be a superbialgebra and let
$$
R_{\sB} = K_0(\com_{\sB})
$$
denote the Grothendieck group of the category $\com_{\sB}$. Thus,
for each $V$ in $\com_{\sB}$, there is an element $[V] \in R_{\sB}$
and each short exact sequence $0 \to U \to V \to W \to 0$ in
$\com_{\sB}$ gives rise to an equation $[V] = [U] + [W]$ in
$R_{\sB}$. The group $R_{\sB}$ is in fact a ring with multiplication
given by the tensor product of $\sB$-comodules. If $\sB$ is
supercommutative as a superalgebra then the ring $R_{\sB}$ is
commutative; see \S\ref{SS:comod}.

Both the ordinary dimension and the superdimension are additive on
short exact sequences and multiplicative on tensor products. Hence
they yield ring homomorphisms
\begin{equation*}
\dim, \sdim \colon R_\sB \to \Z
\end{equation*}
Parts (b) and (c) of Lemma~\ref{L:sch} and formula \eqref{E:sch3a}
have the following immediate consequence:

\begin{cor} \label{C:grothendieck}
The map $[V] \mapsto \schi_V$ yields a well-defined ring
homomorphism $\schi \colon R_{\sB} \to \sB_{\bar 0}$. Furthermore,
the following diagram commutes
$$
\xymatrix{ R_{\sB} \ar[r]^-{\schi} \ar[d]_{\sdim} & \sB_{\bar 0}
\ar[d]^{\e}
\\
\Z \ar[r]_-{\text{\rm can.}} & \k
}
$$
\end{cor}

Forgetting the $\Z_2$-grading,  the corollary also gives the more
familiar version with $\chi$ and $\dim$ in place of $\schi$ and
$\sdim$, respectively.

\subsection{General linear supergroup and Berezinian} \label{SS:superGL}

Let $V$ in $\Vect_{\k}^s$ be finite-dimensional and fix a standard
basis $x_1,\dots,x_d$ with $\dg{i}= \bar 0 \ (i \le p)$ and $\dg{i}
= \bar 1 \ (i>p)$.

\subsubsection{} \label{SSS:EV}

For each supercommutative $\k$-superalgebra $\sR$ we denote by
$\Endo(V)(\sR)$ the set of all $\sR$-linear maps $V \otimes \sR \to
V\otimes \sR$  in $\Vect_{\k}^s$. Using the identification
$\End_{\sR}(V \otimes \sR) \cong \Hom_\k(V,V\otimes \sR) \cong
\End_\k(V) \otimes \sR$ (see \eqref{E:hom}), we may view
$\Endo(V)(\sR)$ as the even subspace of $\End_\k(V) \otimes \sR$:
\begin{equation*}
\Endo(V)(\sR) = \left( \End_\k(V) \otimes \sR \right)_{\bar 0}
\end{equation*}
This defines a functor $\Endo(V)$ from the category of
supercommutative $\k$-superalgebras to the category of semigroups.

\subsubsection{}

Tensoring the supertrace $\str \colon \End_\k(V) \to \k$ of
\eqref{E:str} with $\Id_\sR$, we obtain an $\sR$-linear supertrace
map $\str \colon \End_\k(V) \otimes \sR \to \sR$ in $\Vect_{\k}^s$
which restricts to a map
$\Endo(V)(\sR) \to \sR_{\bar 0}$.
The given standard basis $x_1,\dots,x_d$ of $V$ is an $\sR$-basis of
$V \otimes \sR$. In terms of this basis, an element $\phi\in
\Endo(V)(\sR)$ is given by
\begin{equation} \label{E:Phi}
\phi(x_j)= \sum_{i=1}^d x_i \Phi^i_j \quad \text{with} \quad
\Phi^i_j \in \sR_{\dg{i}+\dg{j}}
\end{equation}
Thus $\phi$ is described by a \emph{supermatrix} $\Phi =
\left(\Phi^i_j\right)$ in \emph{standard form} over $\sR$:
\begin{equation} \label{E:standard}
\Phi= \begin{pmatrix} A & B \\ C & D \end{pmatrix}
\end{equation}
where $A = \left(\Phi^i_j\right)_{i,j\le p}$ and  $D =
\left(\Phi^i_j\right)_{i,j> p}$ are square matrices with entries in
$\sR_{\bar 0}$ while $C, D$ are matrices over $\sR_{\bar 1}$. The
supertrace of $\phi$ is given by
\begin{equation*} \label{E:strPhi}
\str(\phi) = \sum_i (-1)^{\dg{i}} \Phi^i_i = \tr(A) - \tr(D) =:
\str(\Phi)
\end{equation*}

\subsubsection{} \label{SSS:EV2}

The functor $\Endo(V)$ is represented by a supercommutative
$\k$-superbialgebra  which coacts on $V$; this algebra will be
denoted by
\begin{equation*}
\sB = \sO(\Endo(V))
\end{equation*}
Thus, there is a natural isomorphism of $\Endo(V)$ with the functor
$\Hom(\sB,?)$ of parity preserving algebra homorphisms. In
particular, the identity map on $\sB$ corresponds to an element $\xi
\in \Endo(V)(\sB)$. Let $X = (x^i_j)_{d \times d}$ be the matrix of
$\xi$, as in \eqref{E:Phi}. The elements $x^i_j$ have parity $\dg{i}
+ \dg{j}$ and they form a set of supercommuting algebraically
independent generators of $\sB$. In fact, $\sB$ is isomorphic to the
symmetric superalgebra $\Sym(V^* \otimes V)$, with $x^i_j \mapsto
x^i \otimes x_j$, where $\{ x^i \} \subseteq V^*$ is the dual basis
for the given basis of $V$.

We can think of $X$ as the \emph{generic supermatrix} with respect
to the given basis of $V$: any supermatrix $\Phi =
\left(\Phi^i_j\right)$ as in \eqref{E:Phi} comes from an algebra map
$\sB \to \sR$ via $x^i_j \mapsto \Phi^i_j$. The canonical coaction
$\delta \colon V \to V \otimes \sB$, the comultiplication $\Delta$
and the counit $\e$ of $\sB$ are given by
\begin{equation} \label{E:Bformulas}
\begin{aligned}
\delta(x_j) &= \sum_i x_i\otimes x^i_j \\
\Delta(x^i_j) &= \sum_k x^i_k \otimes x^k_j \\
\e(x^i_j) &= \delta^i_j
\end{aligned}
\end{equation}
These formulas can also be written as $\delta(x_1,\dots,x_d) =
(x_1,\dots,x_d)\otimes X$, $\Delta(X) = X \otimes X$ and $\e(X) =
1$.

\subsubsection{} \label{SSS:GLV}

Similarly, $\GL(V)(\sR)$ is defined, for any supercommutative
$\k$-superalgebra $\sR$, as the set of all \emph{invertible}
$\sR$-linear endomorphism of $V\otimes \sR$ in $\Vect_{\k}^s$. The
condition for a supermatrix $\Phi$ in standard form (as in
\eqref{E:standard}) to be invertible is that $A$ and $D$ are
invertible as ordinary matrices over $\sR_{\bar 0}$. In this case,
the inverse of $\Phi$ is given by
\begin{equation*}
\Phi^{-1} =
\begin{pmatrix} (A-BD^{-1}C)^{-1} & -A^{-1}B(D-CA^{-1}B)^{-1}\\
-D^{-1}C(A-BD^{-1}C)^{-1} & (D-CA^{-1}B)^{-1} \end{pmatrix}
\end{equation*}
See Berezin \cite[Theorem 3.1 and Lemma 3.2]{fB87}. The element
\begin{equation} \label{E:ber}
\ber(\Phi):= \det(A)\det(D-CA^{-1}B)^{-1} =
\det(D)^{-1}\det(A-BD^{-1}C)
\end{equation}
is called the \emph{superdeterminant} or \emph{Berezinian} of
$\Phi$; it is an invertible element of $\sR_{\bar 0}$.

The functor $\GL(V)$ is represented by a supercommutative Hopf
superalgebra $\sO(\GL(V))$ which is generated over $\sB =
\sO(\Endo(V))$ by $\det(X_{11})^{-1}$ and $\det(X_{22})^{-1}$, where
$X_{11} = \left( x^i_j \right)_{i,j \le p}$ and $X_{22} = \left(
x^i_j \right)_{i,j > p}$ are the even blocks of the generic
supermatrix $X$. By \cite[Theorem 3.3]{fB87}, the Berezinian
$\ber(X)$ is a group-like element in $\sO(\GL(V))$.

\subsection{Supersymmetric functions and exterior powers}
\label{SS:supersymmetric}

Throughout this section, $V$ will denote a finite-dimensional vector
superspace over $\k$. We assume that the characteristic of $\k$ is
zero.

\subsubsection{} \label{SSS:wedge}
Let
\begin{equation*} \label{E:Snantisymm}
Y_n = \tfrac{1}{n!}\sum_{\sigma \in\Sy_n}\sign(\si) \si \in
\k[\Sy_n]
\end{equation*}
be the \emph{antisymmetrizer} idempotent of the group algebra
$\k[\Sy_n]$ and define
\begin{equation} \label{E:wedge}
\Gr^n V : = \im c_{Y_n} \subseteq V^{\otimes n}
\end{equation}
where $c \colon \k[\Sy_n] \to \End_{\Vect_{\k}^s}(V^{\otimes n})$ is
as in \eqref{E:rep}. Thus, $\Gr^n V$ is the space of antisymmetric
$n$-tensors,
$$
\Gr^n V = \{ y \in V^{\otimes n} \mid c_\si(y) = \sign(\si)y \text{
for all $\si \in \Sy_n$}\}
$$

For later use, we describe an explicit basis of $\Gr^n V$. To this
end, fix a standard basis $x_1,\dots,x_d$ of $V$, with $\dg{i} =
\bar 0$ for $i \le p$ and $\dg{i} = \bar 1$ for $i > p$. Then the
products $x_{\bi} = x_{i_1}\otimes x_{i_2}\otimes\dots\otimes
x_{i_n}$ for sequences $\bi = (i_1,i_2,\dots,i_n) \in
\{1,2,\dots,d\}^n$ form a graded basis of $V^{\otimes n}$ that is
permuted up to a $\pm$-sign by the action of $\Sy_n$ on $V^{\otimes
n}$; see formula \eqref{E:c}:
\begin{equation} \label{E:c'}
c_{\si}(x_{\bi}) = \sign_{\bi}(\si) x_{\si(\bi)}
\end{equation}
with
$$
\sign_{\bi}(\si) =  (-1)^{\sum_{(p,q) \in \inv(\si)}
\dg{i_p}\dg{i_q}} \quad \text{and} \quad \si(\bi) =
(i_{\si^{-1}(1)},i_{\si^{-1}(2)},\dots,i_{\si^{-1}(n)})
$$
Therefore, by elementary properties of monomial group
representations, a $\k$-basis of $\Gr^n V$ is given by the nonzero
elements $c_{Y_n}(x_{\bi})$ where $\bi$ ranges over a transversal
for the $\Sy_n$-action on $\{1,2,\dots,d\}^{n}$. Such a transversal
is provided by the weakly increasing sequences $\bi \in
\{1,2,\dots,d\}^n$. Moreover, for a weakly increasing $\bi$, it is
easily seen from \eqref{E:c'} that $c_{Y_n}(x_{\bi}) = 0$ holds
precisely if $i_\ell = i_{\ell + 1} \le p$ for some $\ell$.
Therefore, a basis of $\Gr^n V$ is given by the elements
$c_{Y_n}(x_{\bi})$ with $\bi = (i_1<i_2<\dots<i_m < i_{m+1} \le
\dots\le i_n) \in \{1,2,\dots,d\}^n$ and $i_m \le p < i_{m+1}$.

In particular,
\begin{equation} \label{E:wedgedim}
\dim_\k \Gr^n V  = \sum_{m + m' = n} \binom{p}{m} \binom{q+m'-1}{m'}
\end{equation}
where $p = \dim_\k V_{\bar 0}$ and $q = \dim_\k V_{\bar 1}$.
Equivalently, the generating power series in $\Z\llbracket t
\rrbracket$ for the sequence $\dim_\k \Gr^n V$ is given by
\begin{equation} \label{E:wedgeseries}
\sum_{n \ge 0} \dim_\k \Gr^n V t^n = \frac{(1+t)^p}{(1-t)^q}
\end{equation}
When $q
> 0$ then all $\Gr^n V$ are nonzero. For additional details on
exterior powers, see, e.g., \cite[Sections I.5 and I.7]{gT04}.

\subsubsection{} \label{SSS:supersymm}

Consider the super bialgebra $\sB = \sO(\Endo(V))$ as defined in
\S\ref{SSS:EV2} and recall that $V$ is in $\com_{\sB}$. The
representation $c \colon \k[\Sy_n] \to
\End_{\Vect_{\k}^s}(V^{\otimes n})$ of \eqref{E:rep} actually has
image in $\End_{\com_{\sB}}(V^{\otimes n})$, since $\sB$ is
supercommutative. Therefore, $\Gr^n V$ also belongs to $\com_{\sB}$
and we can define the \emph{$n^\text{th}$ elementary supersymmetric
function} by
\begin{equation*} \label{E:en}
e_n:= \schi_{\Gr^n V } = \schi(c_{Y_n}) \in \sB_{\bar 0}
\end{equation*}
Here, the equality $\schi_{\Gr^n V } = \schi(c_{Y_n})$ holds by
Lemma~\ref{L:sch}(c).

Similarly, one defines the \emph{$n^\text{th}$ super power sum} by
\begin{equation*} 
p_n:= \schi(c_{(1,2,\dots,n)}) \in \sB_{\bar 0}
\end{equation*}
where $(1,2,\ldots,n) \in \Sy_n$ the cyclic permutation mapping
$1\mapsto 2\mapsto 3\mapsto\ldots\mapsto n\mapsto 1$. In terms of
the generic supermatrix $X$ from \S\ref{SSS:EV2}, one has
$$
p_n = \str(X^n)
$$

Modulo the space spanned by the Lie commutators $fg - gf$ with $f,g
\in \k[\Sy_n]$, the following relation is easily seen to hold in
$\k[\Sy_n]$:
$$
n Y_n \equiv \sum_{i=1}^{n} (-1)^{i-1} (1,2,\dots,i) Y_{n-i}
$$
(with $Y_0=1$). Applying the function $\schi \circ c \colon
\k[\Sy_n] \to \sB_{\bar 0}$ to this relation and using
Lemma~\ref{L:sch}(a),(b), one obtains the \emph{Newton relations}:
\begin{equation*} 
ne_n = \sum_{i=1}^{n}(-1)^{i-1}p_i e_{n-i}
\end{equation*}
Let $t$ be a formal parameter (of parity $\bar 0$) and consider the
generating functions $P(t)= \sum_{n\ge 1}p_n  t^{n-1}$ and $E(t) =
\sum_{n \ge 0} e_n t^n$ in $\sB_{\bar 0}\llbracket t \rrbracket$.
The Newton relations can be written in the form $P(-t)=
\frac{d}{dt}\log E(t)$; see, e.g., \cite[p.~23]{igMac}. Combining
this with the identity
$$
\ber(\exp(tX))=\exp(\str(tX))
$$
due to Berezin (\cite[Chapter~3]{fB87} or \cite[p.~167]{yM97}), one
obtains the following expansion for the characteristic function
$\ber(1 + tX)$ of generic supermatrix $X$:

\begin{prop} \label{P:charfctn}
$\ber(1+ tX) = \sum_{n \ge 0} e_n t^n$
\end{prop}

This proposition is known; see, e.g., Khudaverdian and Voronov
\cite[Prop.~1]{hKtV05}.


\section{Homogeneous superalgebras} \label{S:homog}

\subsection{$N$-homogeneous superalgebras} \label{SS:Nsuperalg}
Let $N$ be an integer with $N \ge 2$. A \emph{homogeneous
superalgebra of degree $N$} or \emph{$N$-homogeneous superalgebra}
is an algebra $\sA$ of the form \eqref{E:superisom} with $V$
finite-dimensional and $R \subseteq V^{\otimes N}$:
\begin{equation*} \label{E:Nsuperalg}
\sA = A(V,R) \cong \T(V)/(R)
\end{equation*}
The assumption $R \subseteq V^{\otimes N}$ implies that, besides the
usual $\Z_2$-grading (``parity''), $\sA$ also has a connected
$\Z_+$-grading (``degree''),
$$
\sA = \bigoplus_{n\ge 0} \sA_n
$$
The algebra $\sA$ is generated by $\sA_1 = V$ and all homogeneous
components $\sA_n$ are finite-dimensional objects of $\Vect_{\k}^s$.
In fact,
\begin{equation} \label{E:Rn}
\sA_n \cong V^{\otimes n}/R_n\quad \text{with}\quad R_n:= (R) \cap
V^{\otimes n} = \sum_{i+j+N=n} V^{\otimes i} \otimes R \otimes
V^{\otimes j}
\end{equation}
Note that $R_n = 0$ for $n < N$; so $\sA_n \cong V^{\otimes n}$ if
$n < N$.

Morphisms of $N$-homogeneous superalgebras $f \colon \sA = A(V,R)
\to \sA' = A(V',R')$  are morphism of superalgebras which also
respect the $\Z_+$-grading. Equivalently, by restricting to degree
$1$, we have a morphism $f_1 \colon \sA_1=V \to \sA'_1=V'$ in
$\Vect_{\k}^s$ whose $N^\text{th}$ tensor power satisfies
$f_1^{\otimes N}(R) \subseteq R'$. Thus, one has a category $\HN$ of
$N$-homogeneous $\k$-superalgebras. Finally, $N$-homogeneous
superalgebras with $N=2$ are called \emph{quadratic} superalgebras;
for $N=3$, they are called \emph{cubic}, etc..

\subsection{Some examples} \label{SS:examples}
In order to explicitly describe a certain $N$-homoge\-neous
superalgebra $\sA = A(V,R)$, we will usually fix a $\Z_2$-graded
$\k$-basis $x_1,\dots,x_d$ of $V = \sA_1$ and denote the the parity
of $x_i$ by $\dg{i}$, as in \S\ref{SS:superspace}. The $x_i$ form a
set of algebra generators for $A$. Following Manin
\cite{yM89},\cite{yM91}, the $d$-tuple $\f = (\dg{1},\ldots,\dg{d})
\in \Z_2^d$ is called the \emph{format} of the basis $\{ x_i \}$.

\begin{example}[Quantum superspace \cite{yM91}] \label{EX:quantum}
For a fixed family $\q$ of scalars $0 \neq q_{ij} \in \k$ $(1 \le i
< j \le d)$ and a given format $\f =(\dg{1},\ldots,\dg{d}) \in
\Z_2^d$ of the basis $x_1,\dots,x_d$, the quadratic superalgebra
$\sA = \Sym_{\q}^{\f}$  is defined as the factor of $\T(V)$
modulo the ideal generated by the elements
\begin{equation} \label{E:affinerels1}
r_i:= x_i \otimes x_i \in (V^{\otimes 2})_{\bar 0} \qquad (\dg{i} =
\bar 1)
\end{equation}
\begin{equation} \label{E:affinerels2}
r_{ij}:= x_j\otimes x_i - q_{ij} (-1)^{\dg{i}\dg{j}} x_i\otimes x_j
\in (V^{\otimes 2})_{\dg{i} + \dg{j}} \qquad(i < j)
\end{equation}
Thus, the algebra $\Sym_\q^{\f}$ is generated by $x_1,\dots,x_d$
subject to the defining relations
\begin{equation*} \label{E:affinerel1}
x_i x_i = 0 \qquad (\dg{i} = \bar 1)
\end{equation*}
and
\begin{equation*} \label{E:affinerel2}
x_j x_i = q_{ij} (-1)^{\dg{i}\dg{j}} x_i x_j \qquad (i < j).
\end{equation*}
In the special case where all $q_{ij}=1$, the algebra
$\Sym_{\q}^{\f}$ is the symmetric superalgebra $\Sym(V)$ of $V$ as
in Example~\ref{EX:SV}.

The ordered monomials of the form $x_1^{m_1}x_2^{m_2}\dots
x_d^{m_d}$, with $\sum_i m_i = n$, $m_i \ge 0$ for all $i$ and $m_i
\le 1$ if $\dg{i} = \bar 1$, form a $\k$-basis of the $n^\text{th}$
homogeneous component of $\Sym_\q^{\f}$. Therefore,
\begin{equation} \label{E:affinedim}
\dim_\k ( \Sym_\q^{\f} )_n = \sum_{r+s=n} \binom{r+p-1}{p-1}
\binom{q}{s}
\end{equation}
where $\dim V_{\bar 0} = p$ and $\dim V_{\bar 1} = q$ as usual.Thus,
the generating series of the dimensions is
\begin{equation*}
\sum_{n\ge 0} \dim_\k ( \Sym_\q^{\f} )_n t^n =
\frac{(1+t)^q}{(1-t)^p}
\end{equation*}
\end{example}

\begin{example}[Yang-Mills algebras \cite{aCmDV03},\cite{aCmDV02}] \label{EX:YM}
Fix a collection of elements $x_1,\dots,x_d$ $(d \ge 2)$, numbered
so as to have parity $\dg{i} = \bar 0$ for $i \le p$ and $\dg{i} =
\bar 1$ for $i > p$. Let $G = (g_{ij}) \in \GL_d(\k)$ be an
invertible symmetric $d \times d$-matrix satisfying $g_{ij} = 0$ if
$\dg{i} \neq \dg{j}$ and consider the cubic superalgebra $\sA$ that
is generated by elements $x_1,\dots,x_d$ subject to the relations
\begin{equation} \label{E:YM1}
\sum_{i,j} g_{ij} [x_i,[x_j,x_k]] = 0 \qquad (k = 1,\dots,d)
\end{equation}
Here $[\,.\,,.\,]$ is the supercommutator \eqref{E:supercomm}. The
algebra $\sA$ will be denoted by $\YM^{p|q}$  $(q = d-p)$. In
particular, the pure even algebra $\YM^{d|0}$ is the ordinary
Yang-Mills algebra introduced in \cite{aCmDV02} while $\YM^{0|d}$ is
the super Yang-Mills algebra as in \cite{aCmDV03}.

As usual, put $V = \sum_i \k x_i$ and let $[\,.\,,.\,]_{\otimes}$
denote the supercommutator in $\T(V)$. Furthermore, put $r_k =
\sum_{i,j} g_{ij} [x_i,[x_j,x_k]_{\otimes}]_{\otimes}$ and $R =
\sum_k \k r_k \subseteq V^{\otimes 3}$; so $\YM^{p|q} = \T(V)/(R)$.
Using the symmetry of $G$, we may replace the $r_k$ by simpler
relations as follows. Choose an invertible $d \times d$-matrix $C =
(c_{ij})$ with $c_{ij} = 0$ if $\dg{i} \neq \dg{j}$ and such that
$C^{\rm tr}GC$ is diagonal, say $\sum_{i,j} c_{ir}g_{ij}c_{js} = g_s
\delta^r_s$. Replace the bases $\{x_i\}$ of $V$ and $\{r_k\}$ of $R$
by the new bases $y_i = \sum_j c^{ij}x_j$ and $s_k = \sum_k c^{k
\ell} r_\ell$ where $C^{-1} = (c^{ij})$. Note that $y_i$ has parity
$\dg{i}$ and $s_k$ has parity $\dg{k}$, the parity of $r_k$. A
simple calculation shows that
$s_k = \sum_{i\neq k} g_i [y_i,[y_i,y_k]_\otimes]_\otimes$.
Thus we obtain the following defining relations for the generators
$y_1,\dots y_d$ of $\YM^{p|q}$:
\begin{equation} \label{E:YM2}
\sum_{i \neq k} g_i [y_i,[y_i,y_k]] = 0 \qquad (k = 1,\dots,d)
\end{equation}

The resulting algebras for $d=2$ are as follows. Putting $x = y_1$
and $y = y_2$ we have two defining relations: $[x,[x,y]]=0$ and
$[y,[y,x]] = 0$. In the pure even case ($\dg{x} = \dg{y} = \bar 0$),
the supercommutators are the ordinary Lie commutators. So
$\YM^{2|0}$ is the enveloping algebra of the Heisenberg Lie algebra;
see \cite[(0.4)]{mAwS87}. In the pure odd case ($\dg{x} = \dg{y} =
\bar 1$), the two relations can be written as $x^2y = y x^2$ and $y
x^2 = x^2 y$. The resulting algebra $\YM^{0|2}$ is a cubic
Artin-Schelter algebra of type $S_1$ \cite[(8.6)]{mAwS87}. Thus,
both unmixed algebras are Artin-Schelter regular of global dimension
$3$. In the mixed case, however ($\dg{x} = \bar 0$, $\dg{y} = \bar
1$), the relations say that $x$ commutes with the Lie commutator
$[x,y]$ while $y$ anticommutes: $y[x,y] = - [x,y]y$. Thus, $[x,y]$
is a normal element of $\YM^{1|1}$ and $\YM^{1|1}/([x,y])$ is a
polynomial algebra in two variables over $\k$. Moreover, the
calculation
\begin{equation*}
[x,y]^2 = [x,[x,y]y] = -[x,y[x,y]] = -[x,y]^2
\end{equation*}
shows that $[x,y]^2 = 0$. Thus, the algebra $\YM^{1|1}$ is
noetherian with Gelfand-Kirillov dimension $2$ and infinite global
dimension.

Returning to the case of general $d \ge 2$, we now concentrate on
the unmixed algebras introduced by Connes and Dubois-Violette. We
will denote these algebras by $\YM^+ = \YM^{d|0}$ and $\YM^- =
\YM^{0|d}$. In all formulas below, $+$ applies to $\YM^+$ and $-$ to
$\YM^-$.
The generators $s_k = \sum_{i\neq k} g_i
[y_i,[y_i,y_k]_\otimes]_\otimes$ of the space of relations $R$ can
be written as $s_k = \sum_\ell y_\ell \otimes m_{\ell k} = \pm
\sum_\ell m_{k\ell} \otimes y_\ell$ with
\begin{equation*}
m_{\ell k} = \begin{cases} g_\ell \left( y_\ell \otimes y_k - (1 \pm
1) y_k \otimes y_\ell \right) \quad &\text{for $\ell \neq k$}\\
\pm \sum_{i \neq k} g_i y_i \otimes y_i &\text{for $\ell = k$}
\end{cases}
\end{equation*}
Thus, putting $Y = (y_1,\dots,y_d)$ and letting $M$ denote the $d
\times d$-matrix over $\YM^\pm$ whose $(\ell,k)$-entry is the image
of $m_{\ell k}$, the defining relations \eqref{E:YM2} can be written
as
\begin{equation} \label{E:YMnew}
Y M  = 0 \qquad \text{or} \qquad  M Y^{\rm tr} = 0
\end{equation}
The defining relations \eqref{E:YM2} for $\sA = \YM^-$ amount to the
even element $\sum_{i} g_{i} y_i^2 \in \sA_2$ being central in
$\sA$.

\end{example}

\begin{example}[$N$-symmetric superalgebra; cf.~\cite{rB01}]
\label{EX:SN}

Let $N \ge 2$ be given and let $V$ be a vector superspace $V$ over a
field $\k$ with $\ch \k = 0$ or $\ch \k > N$. Define
$$
\Sym_N(V) = A(V,R) \quad \text{with } \quad R = \Gr^N V =
c_{Y_N}\left(V^{\otimes N}\right) \subseteq V^{\otimes N}
$$
where $Y_N$ is the antisymmetrizer idempotent of the group algebra
$\k[\Sy_N]$; see \eqref{E:wedge}. This defines a functor
$\Sym_N(\,.\,) \colon \Vect_{\k}^s \to \HN$. Since $2c_{Y_2}$ is the
supercommutator in $\T(V)$, the algebra $\Sym_2(V)$ is just the
symmetric superalgebra $\Sym(V)$ of $V$; see
Example~\ref{EX:quantum}. The algebra $\Sym_N(V)$, for a pure even
space $V = V_{\bar 0}$ and general $N \ge 2$, has been introduced in
\cite{rB01}.

If $2 \le M \le N$ then, viewing $\k[\Sy_M]$ as a subalgebra of
$\k[\Sy_N]$ as usual, the antisymmetrizers of $\k[\Sy_N]$ and
$\k[\Sy_M]$ satisfy $Y_N = Y_M a$ for some $a \in \k[\Sy_N]$.
Therefore,
$$
R = c_{Y_N}\left(V^{\otimes N}\right) \subseteq
c_{Y_M}\left(V^{\otimes N}\right) = c_{Y_M}\left(V^{\otimes
M}\right) \otimes V^{\otimes (N-M)}
$$
This shows that the identity map on $V$ extends to an epimorphism of
superalgebras $\Sym_N(V) \onto \Sym_M(V)$.

Now assume that $\dim_k V = d$ and fix a standard basis
$x_1,\dots,x_d$ of $V$, with $\dg{i} = \bar 0$ for $i \le p$ and
$\dg{i} = \bar 1$ for $i > p$. From the basis for $\Gr^N V$
exhibited in \S\ref{SSS:wedge} we obtain that the algebra
$\Sym_N(V)$ is generated by $x_1,\dots,x_d$ subject to the relations
\begin{equation*} \label{E:SNrels}
\sum_{\si \in \Sy_N} (-1)^{\sum_{(p,q) \in \inv(\si)} 1+
\dg{i_p}\dg{i_q}} x_{i_{\si^{-1}(1)}}x_{i_{\si^{-1}(2)}}\dots
x_{i_{\si^{-1}(N)}} = 0
\end{equation*}
with $1 \le i_1 < i_2 < \dots < i_m \le p = \dim_\k V_{\bar 0} <
i_{m+1} \le \dots \le i_N \le d = \dim_\k V$; see formula
\eqref{E:c'}.

\end{example}

\begin{example} \label{EX:SRN}
The following construction generalizes
Example~\ref{EX:SN}. Fix $N \ge 2$ and $0 \neq q \in \k$ and assume
that condition \eqref{E:Heckess} is satisfied. Given a Hecke
operator $\hR \colon V^{\otimes 2} \to V^{\otimes 2}$ on a vector
superspace $V$
we define the $N$-homogeneous superalgebra
\begin{equation}
\label{E:SRN1}
\Gr_{\hR,N} := A(V,R) \quad \text{with} \quad R = \im
\rho_{\hR}(X_N) \subseteq V^{\otimes N}
\end{equation}
where $X_N \in \He_{N,q}$ is the $q$-symmetrizer \eqref{E:Heckeid1}
and $\rho_\hR$ is the representation \eqref{E:Heckerep} of
$\He_{N,q}$. We also put
\begin{equation}
\label{E:SRN} \Sym_{\hR,N} := \Gr_{-q\hR^{-1},N} = A(V,R) \quad
\text{with} \quad R = \im \rho_{\hR}(Y_N) \subseteq V^{\otimes N}
\end{equation}
where $Y_N \in \He_{N,q}$ is the antisymmetrizer \eqref{E:Heckeid2}.
The algebra $\Sym_N(V)$ in Example~\ref{EX:SN} is identical with
$\Sym_{c_{V,V},N}$ ($q=1$).
\end{example}

\subsection{The dual of a homogeneous superalgebra} \label{SS:dual}
Let $\sA = A(V,R)$ be an $N$-homogeneous superalgebra. The dual
$\sA^!$ of $\sA$ is defined by
$$
\sA^! = A(V^*,R^\perp)
$$
where, $R^\perp \subseteq V^{*\,\otimes N}$ is the (homogeneous)
subspace consisting of all elements that vanish on $R \subseteq
V^{\otimes N}$, using \eqref{E:iso2} in order to evaluate elements
of $V^{*\,\otimes N}$ on $V^{\otimes N}$. Thus, \eqref{E:Rn} takes
the form
\begin{equation} \label{E:Rn!}
\sA^!_n = V^{*\,\otimes n} / R^\perp_n \quad \text{with} \quad
R^\perp_n:= \sum_{i+j+N=n} V^{*\,\otimes i}\otimes R^\perp \otimes
V^{*\,\otimes j}
\end{equation}
Identifying $V^{*\,\otimes n}$ with the linear dual of $V^{\otimes
n}$ via \eqref{E:iso2}, we have  $V^{*\,\otimes i}\otimes R^\perp
\otimes V^{*\,\otimes j} = \left( V^{\otimes j} \otimes R \otimes
V^{\otimes i}\right)^\perp$. Hence,
\begin{equation} \label{E:Rn!2}
R^\perp_n = \left( \bigcap_{i+j+N=n}  V^{\otimes j} \otimes R
\otimes V^{\otimes i} \right)^\perp
\end{equation}

The canonical isomorphism $V \iso V^{**}$ in \eqref{E:**} leads to
an isomorphism $V^{\otimes N} \iso V^{**\,\otimes N}$ which maps $R$
onto $R^{\perp\perp}$. Hence,
\begin{equation} \label{E:!!}
\sA^{!\,!} \cong \sA
\end{equation}
Moreover, if $f \colon \sA = A(V,R) \to \sA' = A(V',R')$ is any
morphism in $\HN$ then the transpose of  $f_1 \colon V \to V'$
induces a morphism $f^! \colon (\sA')^! \to \sA^!$ in $\HN$. Thus,
we have a contravariant quasi-involutive \emph{dualization functor}
$\sA \mapsto \sA^!$, $f \mapsto f^!$ on $\HN$.

\begin{example} \label{EX:unit}
The dual of $A(V,0) = \T(V)$ is $A(V^*,V^{*\,\otimes N})$; so
$$
\T(V)^! = \T(V^*)/\left( V^{*\,\otimes N} \right)
$$
In particular, letting $V = \k$ be the unit object of
$\Vect_{\k}^s$, we have $A(\k,0) = \k[t]$ (polynomial algebra) and
$A(\k,0)^! = \k[d]/(d^N)$, with $t$ and $d$ both having degree $1$
and parity $\bar{0}$.
\end{example}

\begin{example}[Dual of quantum superspace] \label{EX:quantum!}
We will describe the dual $\sA^!$ of quantum superspace $\sA =
\Sym_{\q}^{\f}$; see Example~\ref{EX:quantum}. Fix a homogeneous
$\k$-basis $x_1,\dots,x_d$ with format $\f$ for $V$,
and let $x^1,\dots,x^d$ denote the dual basis of $V^*$; this basis
also has format $\f$. Evaluating an arbitrary element $f =
\sum_{\ell,m} f_{\ell m} x^\ell \otimes x^m \in V^{*\,\otimes 2}$ on
one of the generating relations $r_i,r_{ij}\in R$ in
\eqref{E:affinerels1}, \eqref{E:affinerels2} we obtain $\langle f,
r_{i} \rangle = f_{ii}$ and $\langle f, r_{ij} \rangle = f_{ij} -
q_{ij} (-1)^{\dg{i}\dg{j}}f_{ji}$. Therefore, the space $R^\perp
\subseteq V^{*\,\otimes 2}$ has a basis consisting of the elements
$s^\ell:= x^\ell \otimes x^\ell$ $(\dg{\ell} = \bar 0)$ and
$s^{\ell,k}:= x^\ell\otimes x^k + q_{k \ell} (-1)^{\dg{k}\dg{\ell}}
x^k \otimes x^\ell$ $(k < \ell)$. In summary, $\sA^!$ is generated
by $x^1,\dots,x^d$ subject to the defining relations
\begin{equation*} \label{E:extrel1}
x^\ell x^\ell = 0 \qquad (\dg{\ell} = \bar 0)
\end{equation*}
and
\begin{equation*} \label{E:extrel2}
x^\ell x^k = - q_{k\ell} (-1)^{\dg{k}\dg{\ell}} x^k x^\ell \qquad (k
< \ell).
\end{equation*}
Thus, $\sA^!$ is isomorphic to quantum superspace $\Sym_{\q'}^{\f'}$
with $q_{ij}' = (-1)^{\dg{i}+\dg{j}}q_{ij}$ and $\f' = \f +
(\bar{1},\dots,\bar{1})$ the format obtained from $\f$ by parity
reversal in all components.
\end{example}

\begin{example}[Duals of the Yang-Mills algebras] \label{EX:YM!}
Continuing with the notation of Example~\ref{EX:YM}, we now desribe
the algebra $\sA^!$ for $\sA = \YM^{p|q}$. We assume that $\ch \k =
0$ and work with generators $y_1,\dots,y_d$ of $\sA$ satisfying
\eqref{E:YM2}.

Let $y^1,\dots,y^d$ denote the basis of $V^*$ given by $\langle
y^i,y_j\rangle = \delta^i_j$ and put $\gamma =
\tfrac{1}{d-1}\sum_{i} g_i^{-1} y^i\otimes y^i \in V^{*\otimes 2}$.
Then, for the generators $s_k = \sum_{i\neq k} g_i
[y_i,[y_i,y_k]_\otimes]_\otimes$ of $R$ as in Example~\ref{EX:YM},
one computes
\begin{equation} \label{E:YM!2}
\begin{aligned}
\langle y^a \otimes y^b \otimes y^c, s_k \rangle &=
g_{c}\delta^c_b\delta^a_k +
(-1)^{\dg{b}} g_{b}\delta^b_a\delta^c_k - (-1)^{\dg{a}\dg{k}}(1 + (-1)^{\dg{a}}) g_{a}\delta^a_c\delta^b_k  \\
\langle y^i \otimes \gamma, s_k \rangle &=  \delta^i_k
\end{aligned}
\end{equation}
Therefore, the map $\varphi \mapsto \varphi - \sum_k \langle
\varphi, s_k \rangle y^k \otimes \gamma$ is an epimorphism
$V^{*\otimes 3} \onto R^\perp \subset V^{*\otimes 3}$. We obtain
that the algebra $\sA^!$ is generated by $y^1,\dots,y^d$ subject to
the relations
\begin{equation} \label{E:YM!3}
y^a y^b y^c = (g_{c}\delta^c_b y^a + (-1)^{\dg{b}} g_{b}\delta^b_a
y^c - (-1)^{\dg{a}\dg{b}}(1 + (-1)^{\dg{a}}) g_{a}\delta^a_c y^b)\bg
\end{equation}
where $\bg = \tfrac{1}{d-1}\sum_{i} g_i^{-1} y^i y^i$ is the image
of $\gamma$ in $\sA$.

Since $\sA^!$ is $3$-homogeneous, we clearly have $\sA^!_0 = \k$,
$\sA^!_1 = \bigoplus_i \k y^i = V^*$ and $\sA^!_2 = \bigoplus_{i,j}
\k y^i y^j \cong V^{*\otimes 2}$. By \eqref{E:YM!2}, the elements
$y^a\bg$ form a $\k$-basis of $\sA^!_3 = V^{*\otimes 3}/R^\perp
\cong R^*$. Using the defining relations \eqref{E:YM!3} it is not
hard to see that $\sA^!_4 = \k \bg^2$ and $\sA^!_n = 0$ for $n \ge
5$. If $\sA = \YM^{p|q}$ is of mixed type (i.e., $p \neq 0$ and $q
\neq 0$) then $\bg^2 = 0$.
\end{example}

\begin{example}[Dual of the $N$-symmetric superalgebra]
\label{EX:SN!}%
Recall from Example~\ref{EX:SN} that $\Sym_N(V) = A(V,R)$ with $R =
c_{Y_N}\left(V^{\otimes N}\right)$. Since $Y_N$ is central in
$\k[\Sy_N]$ and stable under the inversion involution $^*$ of
$\k[\Sy_N]$, it follows from \eqref{E:contra2} that
$$
\langle x, c_{Y_N}(y) \rangle = \langle c_{Y_N}(x),y \rangle
$$
holds for all $x \in V^{*\otimes N}$ and $y \in V^{\otimes N}$.
Therefore,
$$
R^\perp = \Ker_{V^{* \otimes N}}(c_{Y_N}) = (1 - c_{Y_N})\left( V^{*
\otimes N} \right)
$$
and so
$$
\Sym_N(V)^! = A\left(V^*,(1 - c_{Y_N})( V^{* \otimes N}) \right)
$$
Note that
\begin{equation} \label{E:SN!}
\bigcap_{i+j+N=n} V^{\otimes i} \otimes R \otimes V^{\otimes j} =
c_{Y_n}\left( V^{\otimes n} \right)
\end{equation}
holds for all $n \ge N$. This follows from \eqref{E:imXn}.
Alternatively, as has been noted in Example~\ref{EX:SN}, we have
$c_{Y_n}\left( V^{\otimes n} \right) \subseteq R \otimes V^{\otimes
(n-N)}$. In the same way, one sees that $c_{Y_n}\left( V^{\otimes n}
\right) \subseteq V^{\otimes i} \otimes R \otimes V^{\otimes j}$
whenever $i+j+N=n$. For the reverse inclusion, note that each $x \in
V^{\otimes i} \otimes R \otimes V^{\otimes j}$ satisfies
$c_{\si_\ell}(x) = -x$ for all transpositions $\si_\ell = (\ell,\ell
+ 1) \in \Sy_n$ with $i < \ell < i+N$. Hence, the left hand side of
\eqref{E:SN!} is contained in
the space of antisymmetric $n$-tensors, $\Gr^n V = c_{Y_n}\left(
V^{\otimes n} \right)$, thereby proving \eqref{E:SN!}. We deduce
from \eqref{E:Rn!}, \eqref{E:Rn!2} and \eqref{E:wedgedim} that
\begin{equation} \label{E:SN!dim}
\dim_\k \Sym_N(V)^!_n = \begin{cases} d^n & \text{if $n < N$}
\\ \sum_{r+s = n} \binom{p}{r} \binom{q+s-1}{s} & \text{if $n \ge N$}
\end{cases}
\end{equation}
where $d = \dim_\k V$, $p = \dim_\k V_{\bar 0}$ and $q = \dim_\k
V_{\bar 1}$.
\end{example}

\subsection{The operations $\circ$ and $\bullet$ on $\HN$}
\label{SS:operations}
Let $\sA = A(V,R)$ and $\sA' = A(V',R')$ be $N$-homogeneous
superalgebras. Following \cite{yM88} and \cite{rBmDVmW} we define
the white and black products $\sA\circ \sA'$ and $\sA \bullet \sA'$
by
\begin{equation*}
\begin{split}
\sA\circ \sA' &= A\left(V \otimes V', c_{\pi_N}\left(R \otimes
V'^{\otimes N} +  V^{\otimes N}\otimes R'\right)\right) \\
\sA\bullet \sA' &= A\left(V \otimes V', c_{\pi_N}\left(R \otimes
 R'\right)\right)
\end{split}
\end{equation*}
where $\pi_N \in \Sy_{2N}$ is the inverse of the permutation
$$
(1,2,\dots,2N) \mapsto (1,N+1, 2, N+2, \dots, k, N+k, \dots, N,2N)
$$
Explicitly, $c_{\pi_N} \colon V^{\otimes N} \otimes V'^{\otimes N}
\tto (V \otimes V')^{\otimes N}$ is the morphism in $\Vect_{\k}^s$
that is given by
\begin{equation} \label{E:pi}
c_{\pi_N}\left( v_1 \otimes \dots v_N \otimes v'_1 \otimes \dots
v'_N \right) =  (-1)^{\sum_i \sum_{j > i} \dg{v'_i}\dg{v}_j}(v_1
\otimes v'_1) \otimes \dots (v_N \otimes v'_N)
\end{equation}
Hence, $c_{\pi_N}\left(R \otimes
 R'\right)$ and $c_{\pi_N}\left(R \otimes V'^{\otimes N} +  V^{\otimes
N}\otimes R'\right)$  are homogeneous subspaces of $(V \otimes
V')^{\otimes N}$ and so $\sA\circ \sA'$ and $\sA \bullet \sA'$
belong to $\HN$.

Under the isomorphism $(V'^* \otimes V^*)^{\otimes N} \iso (V
\otimes V')^{*\,\otimes N}$ which comes from \eqref{E:iso2},
the relations $c_{\pi_N}\left( R'^\perp \otimes R^\perp \right)$ of
$\sA'^! \bullet \sA^!$ map onto the relations $\left(
c_{\pi_N}\left(R \otimes V'^{\otimes N} + V^{\otimes N}\otimes
R'\right) \right)^\perp$ of $(\sA \circ \sA')^!$.
In fact, by \eqref{E:contra2} we have $c^*_{\pi_N} = c_{\pi_N}$,
because $\pi_N\tau = \tau\pi_N$, and so $\langle x, y \rangle =
\langle c_{\pi_N}(x), c_{\pi_N}(y) \rangle$ holds for all $x \in
V'^{*\,\otimes N} \otimes V^{*\,\otimes N}$ and $y \in V^{\otimes N}
\otimes V'^{\otimes N}$. Therefore, canonically,
\begin{equation} \label{E:invol1}
(\sA \circ \sA')^! \cong \sA'^! \bullet \sA^! \quad \text{and} \quad
(\sA \bullet \sA')^! \cong \sA'^! \circ \sA^!
\end{equation}
the two identities being equivalent by \eqref{E:!!}.

By definition of $\circ$, the canonical isomorphisms $\k \otimes V
\cong V \cong V \otimes \k$ in $\Vect_{\k}^s$ give isomorphisms
$A(\k,0) \circ \sA \cong \sA \cong \sA \circ A(\k,0)$ in $\HN$, and
\eqref{E:invol1} yields similar isomorphisms for $\bullet$, with
$A(\k,0)^! = \k[d]/(d^N)$ replacing $A(\k,0) = \k[t]$; see
Example~\ref{EX:unit}.

The supersymmetry isomorphism $c_{V,V'} \colon V \otimes V' \iso V'
\otimes V$ in $\Vect_{\k}^s$ (see \eqref{E:symm}) yields
isomorphisms
\begin{equation} \label{E:symm2}
\sA \circ \sA' \cong \sA' \circ \sA \quad \text{and} \quad \sA
\bullet \sA' \cong \sA' \bullet \sA
\end{equation}
in $\HN$. To see this, note that the following diagram of
isomorphisms in $\Vect_{\k}^s$ commutes:
$$
\xymatrix{ V^{\otimes N} \otimes V'^{\otimes N} \ar[r]^-{c_{\pi_N}}
\ar[d]_{c_{V^{\otimes N},V'^{\otimes N}}} & \left(V \otimes
V'\right)^{\otimes N} \ar[d]^{c_{V,V'}^{\otimes N}}
\\
V'^{\otimes N} \otimes V^{\otimes N} \ar[r]_-{c_{\pi_N}} & \left(V'
\otimes V\right)^{\otimes N} }
$$
with $v_1 \otimes \dots v_N \otimes v'_1 \otimes \dots v'_N \mapsto
(-1)^{\sum_i \sum_{j \ge i} \dg{v'_i}\dg{v}_j}(v'_1 \otimes v_1)
\otimes \dots (v'_N \otimes v_N)$ in both composites. Therefore,
putting $R_{\sA \circ \sA'} = c_{\pi_N}\left(R \otimes V'^{\otimes
N} +  V^{\otimes N}\otimes R'\right)$ and similarly for $R_{\sA'
\circ \sA}$ etc., we have
\begin{equation*}
\begin{split}
c_{V,V'}^{\otimes N}\left( R_{\sA \circ \sA'} \right) &=
\left(c_{\pi_N} \circ c_{V^{\otimes N},V'^{\otimes N}}\right)\left(R
\otimes V'^{\otimes N} +  V^{\otimes N}\otimes
R'\right) \\
&= c_{\pi_N}\left(R' \otimes V^{\otimes N} + V'^{\otimes N}\otimes
R\right) \\
&= R_{\sA' \circ \sA}
\end{split}
\end{equation*}
In the same way, one sees that $c_{V,V'}^{\otimes N}\left( R_{\sA
\bullet \sA'} \right) = R_{\sA' \bullet \sA}$. This proves
\eqref{E:symm2}.

Similarly, the associativity isomorphism $a_{V,V',V''} \colon (V
\otimes V') \otimes V'' \cong V \otimes (V' \otimes V'')$ in
$\Vect_{\k}^s$ leads to isomorphisms
\begin{equation} \label{E:symm3}
(\sA \circ \sA') \circ  \sA'' \cong \sA \circ  (\sA' \circ  \sA'')
\quad \text{and} \quad (\sA \bullet \sA') \bullet  \sA'' \cong \sA
\bullet (\sA' \bullet  \sA'')
\end{equation}
in $\HN$. This is a consequence of the following commutative diagram
of isomorphisms in $\Vect_{\k}^s$: 
$$
\xymatrix{\left( V^{\otimes N} \otimes V'^{\otimes N} \right)
\otimes V''^{\otimes N} \ar[r]_-{c_{\pi_N} \otimes \Id}
\ar[d]_{a_{V^{\otimes N},V'^{\otimes N},V''^{\otimes N}}} & \left(V
\otimes V'\right)^{\otimes N} \otimes V''^{\otimes N}
\ar[r]_-{c_{\pi_N}} & \left((V \otimes V') \otimes
V''\right)^{\otimes N} \ar[d]^{a_{V,V',V''}^{\otimes N}}
\\
V^{\otimes N} \otimes \left( V'^{\otimes N} \otimes V''^{\otimes N}
\right) \ar[r]_-{\Id \otimes c_{\pi_N}} & V^{\otimes N} \otimes (V'
\otimes V'')^{\otimes N} \ar[r]_-{c_{\pi_N}} & \left( V \otimes (V'
\otimes V'')\right)^{\otimes N} }
$$

Finally, the compatibility between the isomorphisms $c_{V,V'}$ and
$a_{V,V',V''}$ (see \S\ref{SS:tensors}) is inherited by the
isomorphisms \eqref{E:symm2} and \eqref{E:symm3} in $\HN$.  To
summarize:

\begin{prop}
The operations $\circ$ and $\bullet$ both make the category $\HN$ of
$N$-homogeneous $\k$-superalgebras into a symmetric tensor category,
with unit objects $A(\k,0) = \k[t]$ for $\circ$ and $A(\k,0)^! =
\k[d]/(d^N)$ for $\bullet$.
\end{prop}

\subsection{The superalgebra map $i \colon \sA \circ \sA' \to \sA \otimes \sA'$}
\label{SS:i}

Let $\sA = A(V,R)$ and $\sA' = A(V',R')$ be objects of $\HN$. The
superalgebra $\sA \otimes \sA'$ is generated by $V \oplus V'$
subject to the relations
$$
R + R' \subseteq (V \oplus V')^{\otimes N} \qquad \text{and} \qquad
[V,V']_\otimes \subseteq (V \oplus V')^{\otimes 2}
$$
where $[\,.\,,.\,]_{\otimes}$ is the supercommutator
\eqref{E:supercomm} in the tensor algebra, as usual. Thus, $\sA
\otimes \sA'$ is not $N$-homogeneous when $N \ge 3$. Nonetheless,
there always is an injective superalgebra homomorphism $i \colon \sA
\circ \sA' \to \sA \otimes \sA'$ which is defined as follows. The
linear embedding $V \otimes V' \into \T(V) \otimes \T(V')$ extends
uniquely to a superalgebra map
\begin{equation} \label{E:iota}
\tilde\iota
\colon \T(V \otimes V') \to \T(V) \otimes \T(V')
\end{equation}
which doubles degrees: the restriction of $\tilde\iota$ to degree
$n$ is the embedding
$$
\T(V \otimes V')_n = (V \otimes V')^{\otimes n}
\stackrel{c_{\pi_n}^{-1}}{\tto} V^{\otimes n} \otimes V'^{\otimes n}
\subseteq (\T(V) \otimes \T(V'))_{2n}
$$
in $\Vect_{\k}^s$, where $c_{\pi_n}$ is as in \eqref{E:pi}. Thus,
$\tilde\iota$ identifies the superalgebra $\T(V \otimes V')$ with
the (super) \emph{Segre product} $\bigoplus_{n\ge 0} V^{\otimes n}
\otimes V'^{\otimes n}$ of $\T(V)$ and $\T(V')$.

The map $\tilde\iota$ sends $R_{\sA \circ \sA'} = c_{\pi_N}\left(R
\otimes V'^{\otimes N} +  V^{\otimes N}\otimes R'\right) \subseteq
(V \otimes V')^{\otimes N}$ to $R \otimes V'^{\otimes N} +
V^{\otimes N}\otimes R'$, the kernel of the canonical epimorphism
$V^{\otimes N} \otimes V'^{\otimes N} \onto \sA_N \otimes \sA'_N$.
Thus:

\begin{prop} \label{P:i}
The algebra map $\tilde\iota$ in \eqref{E:iota} passes down to yield
an injective homomorphism $\k$-superalgebras $i \colon \sA \circ
\sA' \mono \sA \otimes \sA'$ which doubles degree. The image of $i$
is the super Segre product $\bigoplus_{n\ge 0} \sA_n \otimes \sA'_n$
of $\sA$ and $\sA'$.
\end{prop}

\subsection{Internal $\ihom$} \label{SS:ihom}

The isomorphisms \eqref{E:iso1} and \eqref{E:iso2} together with
associativity lead to a functorial isomorphism
\begin{equation*}\label{E:iso3}
\Hom_{\k}(U \otimes V, W^*) \cong \Hom_{\k}(U, (V \otimes W)^*)
\end{equation*}
in $\Vect_{\k}^s$. Explicitly, if  $g \in \Hom_{\k}(U \otimes V,
W^*)$ and $g' \in \Hom_{\k}(U, (V \otimes W)^*)$ correspond to each
other under the above isomorphism then
\begin{equation}\label{E:iso3a}
\langle g(u \otimes v),w\rangle =  \langle g'(u),v \otimes w\rangle
\end{equation}
holds for all $u \in U$, $v \in V$ and $w \in W$.

In particular, by restricting to
$\bar{0}$-components, we have a $\k$-linear isomorphism
\begin{equation}\label{E:iso4}
\Hom_{\Vect_{\k}^s}(U \otimes V, W^*) \cong \Hom_{\Vect_{\k}^s}(U,
(V \otimes W)^*)
\end{equation}
This isomorphism leads to

\begin{prop} \label{P:ihom}
There is a functorial isomorphism
$$
\Hom_{\HN}(\sA \bullet \sB, \sC)
\cong \Hom_{\HN}(\sA, \sC \circ \sB^!)
$$
\end{prop}

\begin{proof}
We follow Manin \cite[4.2]{yM88}. Let $\sA = A(U,R)$, $\sB = A(V,S)$
and $\sC = A(W,T)$ be $N$-homogeneous superalgebras. We will prove
the proposition in the following equivalent form; see \eqref{E:!!} and
\eqref{E:invol1}:
$$
\Hom_{\HN}(\sA \bullet \sB, \sC^!) \cong
\Hom_{\HN}(\sA, (\sB \bullet \sC)^!)
$$
Recall that $\sC^! = A(W^*, T^\perp)$ and $(\sB \bullet \sC)^! =
A((V \otimes W)^*, (c_{\pi_N}(S \otimes T))^\perp)$.
Let $g \colon U \otimes V \to W^*$ be a morphism in $\Vect_{\k}^s$
and let $g' \colon U \to (V \otimes W)^*$ be
the morphism in $\Vect_{\k}^s$ that corresponds to $g$ under
\eqref{E:iso4}. We must show that, for homogeneous subspaces $R
\subseteq U^{\otimes N}$, $S \subseteq V^{\otimes N}$ and $T
\subseteq W^{\otimes N}$,
\begin{equation*} \label{E:show1}
g^{\otimes N}\left( c_{\pi_N}( R \otimes S) \right) \subseteq T^\perp
\quad \Leftrightarrow \quad g'^{\otimes N}\left( R \right)
(c_{\pi_N}(S \otimes T))^\perp
\end{equation*}
Identifying $T^{\perp\perp}$ with $T$ as in \S\ref{SS:dual}, the first
inclusion is equivalent to
\begin{equation} \label{E:show2}
\langle  g^{\otimes N}\left( c_{\pi_N}( R \otimes S) \right),
T \rangle = 0
\end{equation}
while the second inclusion states that
\begin{equation} \label{E:show3}
\langle g'^{\otimes N}(R), c_{\pi_N}(S \otimes T) \rangle = 0
\end{equation}
But \eqref{E:iso3a} shows that \eqref{E:show2} and \eqref{E:show3}
are equivalent, which proves the proposition.
\end{proof}

Proposition~\ref{P:ihom} says that the tensor category
$(\HN,\bullet)$ has an internal $\ihom$ which is given by
$$
\ihom(\sA,\sB) = \sB \circ \sA^!
$$
Explicitly, $\ihom(\sA ,\sB)$ is an object of $\HN$ which represents
the functor $(\HN)^{\text{op}} \to \Sets$, $\sX \mapsto
\Hom_{\HN}(\sX \bullet \sA, \sB)$; so there is an isomorphism of
functors
$$
\Hom_{\HN}(? \bullet \sA, \sB) \cong \Hom_{\HN}(?, \ihom(\sA,\sB))
$$
By general properties of $\ihom$ (see \cite[Def.~1.6]{pDjM}), the
morphism $\Id_{\ihom(\sA,\sB)}$ corresponds to a morphism
\begin{equation} \label{E:mu}
\mu \colon \ihom(\sA,\sB) \bullet \sA \to \sB
\end{equation}
in $\HN$ satisfying the following universal property: for any
morphism $f \colon \sX \bullet \sA \to \sB$ in $\HN$ there exists a
unique morphism $g \colon \sX \to \ihom(\sA,\sB)$ such that
the following diagram commutes:
\begin{equation*} \label{E:univ}
\xymatrix{ \sX \bullet \sA \ar[dr]^-{f} \ar[d]_{g \bullet \Id_{\sA}}
& \\
\ihom(\sA,\sB) \bullet \sA \ar[r]_-{\mu} & \sB }
\end{equation*}
In degree $1$, the map $\mu$ is simply $\Id_V \otimes \ev_U \colon V
\otimes U^* \otimes U \tto V \otimes \k = V$.

From $\ihom(\sB,\sC) \bullet \ihom(\sA,\sB) \bullet \sA
\stackrel{\Id \bullet \mu}{\tto} \ihom(\sB,\sC) \bullet \sB
\stackrel{\mu}{\tto} \sC$ one obtains in this way a composition
morphism
\begin{equation} \label{E:m}
m \colon \ihom(\sB,\sC) \bullet \ihom(\sA,\sB) \to \ihom(\sA,\sC)
\end{equation}
in $\HN$. The morphisms $\mu$ and $m$ satisfy the obvious
associativity properties.

\subsection{The superbialgebra $\eA$} \label{SS:eA}

Following Manin \cite[4.2]{yM88} we define
$$
\hom(\sA,\sB) = \ihom(\sA^!,\sB^!)^! = \sA^! \bullet \sB
$$
for $\sA$, $\sB$ in $\HN$. Applying the dualization functor to
\eqref{E:mu}, \eqref{E:m} and recalling \eqref{E:invol1}, we obtain
morphisms
\begin{align*}
\delta_\circ &\colon \sA \to \sB \circ \hom(\sB,\sA)  \\ 
\Delta_\circ &\colon \hom(\sA,\sC) \to \hom(\sA,\sB) \circ
\hom(\sB,\sC) 
\end{align*}
in $\HN$. The associativity properties of $\mu$ and $m$ translate
into corresponding coassociativity properties for $\delta_\circ$ and
$\Delta_\circ$. Following $\delta_\circ$ and $\Delta_\circ$ by the
algebra map $i$ of Proposition~\ref{P:i}, we obtain superalgebra
maps
\begin{align}
\delta &\colon \sA \to \sB \otimes\hom(\sB,\sA) \label{E:delta} \\
\Delta &\colon \hom(\sA,\sC) \to \hom(\sA,\sB) \otimes \hom(\sB,\sC)
\label{E:Delta}
\end{align}

Now take $\sA = \sB = \sC = A(V,R)$ and put $\eA = \hom(\sA,\sA)$;
so
\begin{equation} \label{E:eA}
\eA = \sA^! \bullet \sA = A(V^* \otimes V, c_{\pi_N}(R^\perp \otimes
R))
\end{equation}
Then \eqref{E:Delta} yields a
coassociative superalgebra map
\begin{equation*} \label{E:Delta2}
\Delta \colon \eA \to \eA \otimes \eA
\end{equation*}
Moreover, by Proposition~\ref{P:ihom}, the morphism $\sA^!
\stackrel{\Id}{\tto} \sA^! \cong \k[t] \circ \sA^!$ corresponds to a
morphism $\eA = \sA^! \bullet \sA \to \k[t]$ in $\HN$. Following
this morphism by the map $t \mapsto 1$ we obtain a superalgebra map
\begin{equation*}
\varepsilon \colon \eA \to \k
\end{equation*}
which in degree $1$ is the usual evaluation pairing $\ev_{V} \colon
V^* \otimes V \to \k$ in $\Vect_{\k}^s$. Finally, \eqref{E:delta}
provides us with a superalgebra map
\begin{equation} \label{E:delta2}
\delta_{\sA} \colon \sA \to \sA \otimes \eA
\end{equation}
Note that $\delta_{\sA}$ maps the degree $n$-component of $\sA$
according to
\begin{equation} \label{E:delta3}
\sA_n \stackrel{\delta_\circ}{\tto} (\sA \circ \eA)_n
\stackrel{i}{\tto} \sA_n \otimes (\eA)_n \into \sA_n \otimes \eA
\end{equation}
Fixing a graded $\k$-basis $x_1,\dots,x_d$ of $V$ and denoting the
dual basis of $V^*$ by $x^1,\dots,x^d$ as before, $\eA$ has algebra
generators
\begin{equation} \label{E:eAgens}
z^i_j : = x^i \otimes x_j
\end{equation}
of degree-$1$ and parity $\dg{i} + \dg{j}$. In terms of these
generators, the maps $\varepsilon$, $\delta_{\sA}$ and $\Delta$ are
given by
\begin{equation} \label{E:eAformulas}
\begin{aligned}
\varepsilon(z^j_i) &= \delta^j_i & &\text{or} & \e(Z) &= 1\\
\delta_{\sA}(x_j)  &= \sum_i x_i \otimes z^i_j & &\text{or} &
\delta_{\sA}(x_1,\dots,x_d) &=
(x_1,\dots,x_d) \otimes Z\\
\Delta(z^i_j) &= \sum_k z^i_k \otimes z^k_j & &\text{or} & \Delta(Z)
&= Z \otimes Z
\end{aligned}
\end{equation}
where $Z = ( z^i_j )_{d \times d}$.

\begin{prop} \label{P:eA}
Let $\sA = A(V,R)$ be an $N$-homogeneous $\k$-superalgebra.
\begin{enumerate}
\item
With $\Delta$ as comultiplication and $\varepsilon$ as counit,
the superalgebra $\eA$ becomes a superbialgebra. Moreover,
$\delta_{\sA}$ makes $\sA$ into a graded right $\eA$-comodule
superalgebra.
\item
Given any $\k$-superalgebra $\sB$ and a morphism of superalgebras
$\delta \colon \sA \to \sA \otimes \sB$ satisfying $\delta(V)
\subseteq V \otimes \sB$, there is a unique morphism of
superalgebras $\varphi \colon \eA \to \sB$ such that the following
diagram commutes:
\begin{equation*}
\xymatrix{ \sA \ar[r]^-{\delta} \ar[dr]_{\delta_{\sA}}
& \sA \otimes \sB \\
& \sA \otimes \eA  \ar[u]_-{\Id_{\sA} \otimes \varphi}}
\end{equation*}
\end{enumerate}
\end{prop}

The proposition is proved as in \cite[\S5]{yM88} or \cite[Theorem
3]{rBmDVmW}.

\begin{example} \label{EX:endfree}
When $\sA = A(V,0) = \T(V)$, we have $\eA = A(V^* \otimes V,0) =
\T(V^* \otimes V)$; so
\begin{equation*} \label{E:endTV}
\en \T(V) = \T(V^* \otimes V)
\end{equation*}
the free superalgebra generated by the elements $z^i_j$ in
\eqref{E:eAgens}.
\end{example}

\begin{example} \label{EX:endSN}
By Examples~\ref{EX:SN} and \ref{EX:SN!}, we have
\begin{equation*}
\en \Sym_N(V) = A\left(V^* \otimes V, c_{\pi_N}\left(
(1-c_{Y_N})(V^{*\,\otimes N}) \otimes c_{Y_N}(V^{\otimes N})
\right)\right)
\end{equation*}
For example, the algebra $\en \Sym_2(V)$ is generated by the
elements $z^i_j$ with parity $\dg{i} + \dg{j}$ subject to the
relations
\begin{equation*}
[z^{i_1}_{j_1}, z^{i_2}_{j_2}] + (-1)^{\dg{i_1}\dg{i_2} + (\dg{i_1}
+ \dg{i_2})\dg{j_1}} [z^{i_2}_{j_1}, z^{i_1}_{j_2}] = 0
\end{equation*}
where $[\,.\,,\,.\,]$ is the supercommutator \eqref{E:supercomm}.
This algebra is highly noncommutative, even for a pure even space
$V$.

Let $\sO(\Endo(V)) = \Sym(V^* \otimes V)$ be the supercommutative
superbialgebra as in \S\ref{SSS:EV2}, with generators $x^i_j$. There
is a map of superbialgebras
\begin{equation} \label{E:endSN}
\varphi \colon \en \Sym_N(V) \to \sO(\Endo(V))\,, \qquad z^i_j
\mapsto x^i_j
\end{equation}
Indeed, write $\sB = \sO(\Endo(V))$ for brevity and recall the
coaction $\delta \colon V \to V \otimes \sB$, $ x_j \mapsto \sum_i
x_i \otimes x^i_j$ from \eqref{E:Bformulas}. Since $c_{Y_N} \in
\End_{\com_{\sB}}(V^{\otimes N})$ (see \S\ref{SSS:supersymm}), the
map $\delta$ extends to a map of superalgebras
\begin{equation*}
\delta \colon \Sym_N(V) \to \Sym_N(V) \otimes \sB
\end{equation*}
Therefore, Proposition~\ref{P:eA}(b) yields the desired $\varphi$.
Note that the coaction of $\en \Sym_N(V)$ on $V$, when restricted
along $\varphi$, becomes the canonical coaction of $\sO(\Endo(V))$
on $V$; see \eqref{E:Bformulas} and \eqref{E:eAformulas}.
\end{example}


\section{$N$-Koszul superalgebras} \label{S:koszul}

Throughout this section, we fix an $N$-homogeneous superalgebra $\sA
= A(V,R)$.

\subsection{The graded dual $\sA^{!*}$} \label{SS:gradeddual}

The graded dual
$$
\sA^{!*} = \bigoplus_n \sA^{!\ *}_n
$$
of $\sA^!$ has a natural structure of a graded right $\eA$-comodule.
Indeed, the linear dual $\sA^{!\ *}_{n}$ of the degree $n$-component
of $\sA^!$ embeds into $V^{\otimes n}$ as follows. Recall from
\eqref{E:Rn!2} that
\begin{equation} \label{E:A!*}
\sA^{!\ *}_{n} = \begin{cases} V^{\otimes n} & \text{if $n < N$} \\
\bigcap_{i+j+N=n} V^{\otimes i}\otimes R \otimes V^{\otimes j} &
\text{if $n \ge N$}
\end{cases}
\end{equation}
This identification makes the graded dual $ \sA^{!*}$ into a graded
right $\eA$-comodule. For, by \eqref{E:delta3} the coaction
$\delta_{\sA}$ restricts in degree $1$ to a map $V \to V \otimes
\eA$ which makes $\T(V)$ into a graded right $\eA$-comodule
superalgebra. The structure map $\T(V) \to \T(V) \otimes \eA$ sends
$R \to R \otimes \eA$. Therefore, each $V^{\otimes i}\otimes R
\otimes V^{\otimes j}$ is a $\eA$-subcomodule of $V^{\otimes
(i+j+N)}$, and hence $\sA^{!\ *}_{n}$ is a $\eA$-subcomodule of
$V^{\otimes n}$. Finally, for all $n \ge 0$,
\begin{equation} \label{E:A!*2}
\sA^{!\ *}_{n+1} \subseteq V \otimes \sA^{!\ *}_{n} \quad \text{and}
\quad \sA^{!\ *}_{n+N} \subseteq V^{\otimes N} \otimes \sA^{!\
*}_{n} \cap R \otimes V^{\otimes n} = R \otimes \sA^{!\ *}_{n}
\end{equation}

\subsection{The Koszul complex} \label{SS:KA}

The map
\begin{equation*} \label{E:KA}
\begin{split}
\sA \otimes V^{\otimes (i+1)} &\to  \sA \otimes V^{\otimes
i} \\
a \otimes (v_1 \otimes \dots \otimes v_{i+1}) &\mapsto  av_1 \otimes
(v_2 \otimes \dots \otimes v_{i+1})
\end{split}
\end{equation*}
is a morphism in the category $\comod_{\eA}$ of right
$\eA$-comodules, because the $\eA$-coaction $\delta_{\sA}$ in
\eqref{E:delta2} is a superalgebra map. Furthermore, this map is a
left $\sA$-module map which preserves total degree, and it restricts
to a map of $\eA$-subcomodules
\begin{equation*}
d \colon \sA \otimes \sA^{!\ *}_{i+1} \to \sA V \otimes \sA^{!\
*}_{i} \into \sA \otimes \sA^{!\ *}_{i}
\end{equation*}
which is the $\sA$-linear extension of the embedding \eqref{E:A!*2}.
The map $d^N$ sends $\sA^{!\ *}_{i+N}$ to $\sA R \otimes \sA^{!\
*}_{i} = 0$; so $d^N = 0$. In other words, we have an
\emph{$N$-complex}
\begin{equation} \label{E:KA1}
K(\sA) \colon \ldots \stackrel{d}{\tto} \sA \otimes \sA^{!\ *}_{i+1}
\stackrel{d}{\tto} \sA \otimes \sA^{!\ *}_{i} \stackrel{d}{\tto}
\ldots \stackrel{d}{\tto} \sA \tto 0
\end{equation}
in $\comod_{\eA}$ consisting of graded-free left $\sA$-modules and
$\sA$-module maps which preserve total degree. Therefore, $K(\sA)$
splits into a direct sum of $N$-complexes $K(\sA)^n =
\bigoplus_{i+j=n} \sA_i \otimes \sA^{!\ *}_{j}$ in $\com_{\eA}$.

Following \cite{rBmDVmW}, the \emph{Koszul complex} $\K(\sA)$
defined by Berger in \cite{rB01} can be described as the following
contraction of $K(\sA)$:
\begin{equation} \label{E:KA2}
\K(\sA) \colon \ldots \stackrel{d^{N-1}}{\tto} \sA \otimes \sA^{!\
*}_{N+1} \stackrel{d}{\tto} \sA \otimes \sA^{!\ *}_{N}
\stackrel{d^{N-1}}{\tto} \sA \otimes \sA^{!\ *}_{1}
\stackrel{d}{\tto} \sA \tto 0
\end{equation}
This is an ordinary complex in $\comod_{\eA}$ which splits into a
direct sum of complexes $\K(\sA)^n$ in $\com_{\eA}$. The
$i^\text{th}$ components of $\K(\sA)$ and of $\K(\sA)^n$ are given
by
\begin{equation*} \label{E:jump}
\K(\sA)_i = \sA \otimes \sA^{!\ *}_{\ju(i)}\quad\text{and}\quad
\K(\sA)^n_i = \sA_{n - \ju(i)} \otimes \sA^{!\ *}_{\ju(i)}
\end{equation*}
with $\ju(i) = \ju_N(i)$ as in \eqref{E:jump1}. The differential on
$\K(\sA)$ is
\begin{equation*} \label{E:differential}
\delta_i \colon \K(\sA)_i \to \K(\sA)_{i-1} \quad \text{where} \quad
\delta_i = \begin{cases} \ d^{N-1} & \text{for $i$ even} \\
\ d & \text{for $i$ odd}
\end{cases}
\end{equation*}
Writing $\sA_+ = \bigoplus_{n>0} \sA_n = \sA V$ as usual, we have
\begin{equation*} \label{E:kernel}
\Ker \delta_i \subseteq \sA_+ \K(\sA)_i
\end{equation*}
for all $i$. Indeed, this is clear for odd $i$, since $\delta_i = d$
is injective on $\sA^{!\ *}_{\ju(i)}$. For even $i$, the restriction
of $\delta_i = d^{N-1}$ to $\sA^{!\ *}_{\ju(i)}$ is given by
$d^{N-1} \colon \sA^{!\ *}_{\ju(i)} = \sA^{!\ *}_{\ju(i-1) + N-1}
\into V^{\otimes (N-1)} \otimes \sA^{!\ *}_{\ju(i-1)} \iso \sA_{N-1}
\otimes \sA^{!\ *}_{\ju(i-1)} \into \sA \otimes \sA^{!\
*}_{\ju(i-1)}$ where the first embedding comes from \eqref{E:A!*2}.

Since $\sA^{!\ *}_{\ju(1)} = \sA^{!\ *}_{1} = V$ and $\sA^{!\
*}_{\ju(2)} = \sA^{!\ *}_{N} = R$ by \eqref{E:A!*}, the start of the
Koszul complex, augmented by the canonical map
$$
\sA \onto \k =\sA/\sA_+
$$
is as follows:
\begin{equation} \label{E:start}
\sA \otimes R \stackrel{\delta_2}{\tto} \sA \otimes V
\stackrel{\delta_1 = \text{mult}}{\tto} \sA \tto \k \tto 0
\end{equation}
This piece is easily seen to be exact:
writing $\sA = \T(V)/I$ with $I = (R) = I \otimes V + \T(V) \otimes
R$ as in \eqref{E:I}, the map $\T(V)_+ = \T(V) \otimes V \onto \sA
\otimes V \stackrel{\delta_1}{\onto} \sA_+$ has kernel $I$. Thus,
$\Ker \delta_1 = I/I\otimes V = \im \delta_2$. Hence \eqref{E:start}
is the start of the minimal graded-free resolution of the left
$\sA$-module $\k$.

\subsection{$N$-homogeneous Koszul superalgebras} \label{SS:koszul}

Recall from the Introduction that an $N$-homoge\-neous superalgebra
$\sA$ is called \emph{$N$-Koszul} if $\Tor^{\sA}_i(\k,\k)$ is
concentrated in degree $\ju_N(i)$ for all $i \ge 0$. By
\cite[Proposition 2.12]{rB01} or \cite[Theorem 2.4]{rBnM06}, this
happens exactly if the Koszul complex $\K(\sA)$ is exact in degrees
$i>0$ and in view of \eqref{E:start}, this amounts to exactness of
$\K(\sA)$ in degrees $i \ge 2$. In this case,
\begin{equation*}
\K(\sA) \tto
\k \tto 0
\end{equation*}
is the minimal graded-free resolution of
the trivial left $\sA$-module $\k$.

The Yoneda $\Ext$-algebra $E(\sA) = \bigoplus_{i \ge 0}
\Ext^i_{\sA}(\k,\k)$ of an $N$-Koszul superalgebra $\sA$ has the
following description in terms of the dual algebra $A^!$:
\begin{equation*} \label{E:Ext}
\Ext^i_{\sA}(\k,\k) \cong \sA^!_{\ju(i)} \qquad (i \ge 0)
\end{equation*}
Moreover, identifying $\Ext^i_{\sA}(\k,\k)$ and $\sA^!_{\ju(i)}$,
the Yoneda product $f\cdot g$ and the $\sA^!$-product $fg$ for $f
\in \Ext^i_{\sA}(\k,\k) = \sA^!_{\ju(i)}$ and $g \in
\Ext^j_{\sA}(\k,\k) = \sA^!_{\ju(j)}$ are related by $f \cdot g =
(-1)^{ij}fg$ when $N=2$, and
$$
f \cdot g = \begin{cases} fg \quad&\text{if $i$ or $j$ is even} \\
0 \quad&\text{if $i$ and $j$ are both odd}
\end{cases}
$$
for $N > 2$; see \cite[Theorem 9.1]{eGetal04}, \cite[Proposition
3.1]{rBnM06}.

\begin{example}
Quadratic algebras having a PBW-basis are $2$-Koszul; see, e.g.,
\cite[Chap. 4, Theorem 3.1]{PP05}. This applies in particular to
quantum superspace $\sA = A_{\q}^{\f}$; see
Example~\ref{EX:quantum}. A PBW-basis in this case is given by the
collection of ordered monomials $x_1^{m_1}x_2^{m_2}\dots x_d^{m_d}$
with $m_i \ge 0$ for all $i$ and $m_i \le 1$ if $\dg{i} = \bar 1$,
as in Example~\ref{EX:quantum}. For a more general result, see
\cite[Chap. 4, Theorem 8.1]{PP05}.
\end{example}

\begin{example}
The unmixed Yang-Mills algebras $\sA = \YM^\pm$ (see
Example~\ref{EX:YM}) were shown to be $3$-Koszul in \cite{aCmDV02},
\cite{aCmDV03}. Indeed, letting $\sA[\ell]$ denotes the shift of
$\sA$ that is defined by $\sA[\ell]_n = \sA_{\ell+n}$, the defining
relations for $\sA$ in the form \eqref{E:YMnew} imply that the
following complex of graded-free left $\sA$-modules is exact:
\begin{equation} \label{YMcomplex}
0 \tto \sA[-4] \stackrel{\cdot Y}{\tto} \sA[-3]^d \stackrel{\cdot
M}{\tto} \sA[-1]^d \stackrel{\cdot Y^{\rm tr}}{\tto} \sA \tto \k
\tto 0
\end{equation}
The piece $\sA[-3]^d \stackrel{\cdot M}{\tto} \sA[-1]^d
\stackrel{\cdot Y^{\rm tr}}{\tto} \sA \tto \k \tto 0$ is identical
with \eqref{E:start}. Therefore, \eqref{YMcomplex} is the minimal
graded-free resolution of $\k$. The resolution shows that each
$\Tor^{\sA}_i(\k,\k)$ is concentrated in degree $\ju_3(i)$, and
hence $\sA$ is $3$-Koszul. It also follows that \eqref{YMcomplex} is
isomorphic to $\K(\sA) \to \k \to 0$. In particular,
\eqref{YMcomplex} confirms the dimensions of the corresponding
components $\sA_n^!$ in Example~\ref{EX:YM!}. As has been pointed
out in \cite{aCmDV02}, \cite{aCmDV03}, it follows from
\eqref{YMcomplex} that the Hilbert series $H_{\sA}(t) = \sum_{n \ge
0} \dim_\k \sA_n\,t^n $ of $\sA = \YM^\pm$ has the form
\begin{equation*}
\begin{aligned}
H_{\sA}(t) &= \frac{1}{1 - dt + dt^3 - t^4 }\\
&= \frac{1}{(1 - t^2)(1 - dt + t^2)}
\end{aligned}
\end{equation*}
If $d > 2$ then the series has a pole in the interval $(0,1)$, and
hence $\dim_\k\sA_n$ grows exponentially with $n$. Therefore, $\sA$
is not noetherian in this case; see Stephenson and Zhang
\cite{dStjZ97}.

The mixed Yang-Mills algebras $\sA = \YM^{p|q}$ with $p \neq 0$ and
$q \neq 0$, on the other hand, are never $3$-Koszul. For $\YM^{1|1}$
this follows from the description given in Example~\ref{EX:YM}: this
algebra has infinite global dimension. In general, one can check
that the so-called extra condition (see \eqref{E:ec} below) fails
for $\sA$, and so $\sA$ cannot be Koszul, by \cite[Prop.~2.7]{rB01}.
\end{example}

\begin{example}
It has been shown in \cite[Theorem 3.13]{rB01} that the
$N$-symmetric algebra $\Sym_N(V)$ of a pure even space $V$ over a
field of characteristic $0$ is $N$-Koszul. An extension of this
result will be offered in Theorem~\ref{T:HeckeKoszul} below.
\end{example}

\subsection{Confluence and Koszulity} \label{SS:confluence}

For the convenience of the reader, we recall the notions of
reduction operators and confluence and their relation to the Koszul
property. Complete details can be found in Berger \cite{rB98},
\cite{rB01}.

Let $V$ in $\Vect_{\k}^s$ be given along with a graded basis $X =
\{x_1, \dots x_d\}$ that is ordered by $x_1
> x_2 > \dots > x_d$. The tensors (``monomials") $x_{\bi} =
x_{i_1}\otimes x_{i_2}\otimes\dots\otimes x_{i_N}$ for $\bi =
(i_1,i_2,\dots,i_N) \in \{1,2,\dots,d\}^N$ form a basis of
$V^{\otimes N}$ which will be given the lexicographical ordering. An
\emph{$X$-reduction operator} on $V^{\otimes N}$ is a projection $S
\in \End_{\Vect_{\k}^s}(V^{\otimes N})$ such that either $S(x_\bi) =
x_\bi$ or $S(x_\bi) < x_\bi$ holds for each $\bi$, where the latter
inequality means that $S(x_\bi)$ is a linear combination (possibly
$0$) of monomials $< x_\bi$. The monomials $x_\bi$ satisfying
$S(x_\bi) = x_\bi$ are called \emph{$S$-reduced}, all other
monomials are \emph{$S$-nonreduced}. We denote by $\Red(S)$ and
$\NRed(S)$) the (super) subspaces of $V^{\otimes N}$ that are
generated by the $S$-reduced monomials and the $S$-nonreduced
monomials, respectively; so $V^{\otimes N} = \Red(S) \oplus
\NRed(S)$ and $\im(S) = \Red(S)$.

Let $\mathcal{L}_X(V^{\otimes N})$ denote the collection of all
$X$-reduction operators on $V^{\otimes N}$. The proof of
\cite[Theorem 2.3]{rB98} shows that the application $S \mapsto
\Ker(S)$ is a bijection between $\mathcal{L}_X(V^{\otimes N})$ and
the set of all super subspaces of $V^{\otimes N}$. Hence
$\mathcal{L}_X(V^{\otimes N})$ inherits a lattice structure: for
$S,S' \in \mathcal{L}_X(V^{\otimes N})$ one has $X$-reduction
operators $S \me S'$ and $S \jo S'$ on $V^{\otimes N}$ which are
defined by
\begin{eqnarray*}
\Ker(S \me S') &=& \Ker(S) + \Ker(S')\\
\Ker(S \jo S') &=& \Ker(S) \cap \Ker(S')
\end{eqnarray*}
A pair $(S,S')$ of $X$-reduction operators is said to be
\emph{confluent} if
$$
\Red(S \jo S') = \Red(S) + \Red(S')
$$
Since the inclusion $\supseteq$ is always true, confluence of
$(S,S')$ is equivalent to the inequality
\begin{equation} \label{E:conf}
\dim_\k \im(S \jo S') \le \dim_\k (\im(S) + \im(S'))
\end{equation}

Let $n \ge N$. Any $X$-reduction operator $S$ on $V^{\otimes N}$
gives rise to $X$-reduction operators $S_{n,i}$ on $V^{\otimes n}$
which are defined by
$$
S_{n,i}:= \Id_{V^{\otimes i}} \otimes S \otimes \Id_{V^{\otimes j}}
\qquad (i+j+N = n)
$$
A monomial $x_\bi = x_{i_1}\otimes x_{i_2}\otimes\dots\otimes
x_{i_n}$ of length $n \ge N$ is said to be $S$-reduced if $x_\bi$ is
$S_{n,i}$-reduced for all $i$, that is, if every connected
submonomial of $x_\bi$ of length $N$ is $S$-reduced.

Now let  $\sA = A(V,R)$ be an $N$-homogeneous superalgebra, and let
$S$ be the $X$-reduction operator on $V^{\otimes N}$ such that
$\Ker(S) = R$. The algebra $\sA$ is said to be \emph{$X$-confluent}
if the pairs
$\left(S_{N+i,i},S_{N+i,0}\right)$
of $X$-reduction operators on $V^{\otimes N+i}$ are confluent for
$i=1,\ldots,N-1$. By \eqref{E:conf} this amounts to the inequalities
\begin{equation}\label{E:confluence}
\dim_\k \im(S_{N+i,i} \jo S_{N+i,0}) \le \dim_\k
\left(\im(S_{N+i,i}) + \im(S_{N+i,0})\right)
\end{equation}
being satisfied for $i=1,\ldots,N-1$.

Following Berger \cite{rB01}, we denote by $\scrT_n$ the lattice of
super subspaces of $V^{\otimes n}$ that is generated by the
subspaces
\begin{equation}\label{E:Rni}
R_{n,i}:= V^{\otimes i} \otimes R \otimes V^{\otimes j} =
\Ker(S_{n,i}) \qquad (i+j+N = n)
\end{equation}
The superalgebra $\sA$ is said to be \emph{distributive} if the
lattices $\scrT_n$ are distributive for all $n$, that is, $C \cap
(D+E) = (C \cap D) + (C \cap E)$ holds for all $C,D,E \in \scrT_n$.

The following proposition states the operative facts concerning
Koszulity for our purposes. Part (a) is identical with
\cite[Thm.~3.11]{rB01} while (b) is \cite[Prop.~3.4]{rB01}.

\begin{prop} \label{P:confluence}
Let  $\sA = A(V,R)$ be an $N$-homogeneous superalgebra.
\begin{enumerate}

\item If $\sA$ is $X$-confluent for some totally ordered graded
basis $X$ of $V$ then $\sA$ is distributive. Moreover, letting $S$
denote the $X$-reduction operator on $V^{\otimes N}$ such that
$\Ker(S) = R$, the classes in $\sA$ of the $S$-reduced monomials
$x_{i_1}\otimes x_{i_2}\otimes\dots\otimes x_{i_n}$ with $x_{i_j}
\in X$ form a $\k$-basis of $\sA_n$ for all $n \ge N$.

\item Assume that $\sA$ is distributive and the following ``extra condition" is
satisfied
\begin{equation} \label{E:ec}
R_{n+N,0} \cap R_{n+N,n} \subseteq R_{n+N,n-1} \qquad (2 \le n \le
N-1)
\end{equation}
Then $\sA$ is $N$-Koszul.

\end{enumerate}

\end{prop}

After these preparations, we are now ready to prove the following
result. The quadratic case $N=2$ is due to Gurevich \cite{dG90}; see
also Wambst \cite{mW93}.

\begin{thm} \label{T:HeckeKoszul}
Let $N \ge 2$ and $0 \neq q \in \k$ and assume that $[n]_q \neq 0$
for all $n \ge 1$. Then, for every Hecke operator $\hR$ associated
with $q$, the $N$-homogeneous superalgebra $\Gr_{\hR,N}$ defined in
\eqref{E:SRN1} is $N$-Koszul.
\end{thm}

\begin{proof}
Put $\sA = \Gr_{\hR,N}$ and recall that $\sA = A(V,R)$ with
$$
R = \im \rho_{N,\hR}(X_N) \subseteq V^{\otimes N}
$$
The extra condition \eqref{E:ec} is a consequence of equation
\eqref{E:imXn}. Indeed, \eqref{E:imXn} implies that the spaces
$R_{n,i}$ in \eqref{E:Rni} have the form
\begin{equation}\label{E:imXn2}
R_{n,i} = \bigcap_{s=i+1}^{i+N-1} \im(\rho_{n,\hR}(T_s)+1) \subseteq
V^{\otimes n}
\end{equation}
Applying \eqref{E:imXn2} with $\rho = \rho_{n+N,\hR}$ we see that
the left hand side of \eqref{E:ec} is identical to
$$
\bigcap_{i=1}^{N-1}\im(\rho(T_i) + 1) \cap
\bigcap_{i=n+1}^{n+N-1}\im(\rho(T_i) + 1) =
\bigcap_{i=1}^{n+N-1}\im(\rho(T_i) + 1)
$$
where the equality holds because $n+1 \le N$. The last expression is
clearly contained in $\bigcap_{i=n}^{n+N-2}\im(\rho(T_i) + 1)$,
which is identical to the right hand side of \eqref{E:ec}. This
establishes the extra condition \eqref{E:ec}.

In order to prove the distributivity of $\sA$, we follow the
approach taken in \cite{phH97}. We first prove the claim for the
standard solution $\hR^{DJ}$, i.e., the operator given in
Example~\ref{EX:HeckDJ} with $d=p$ and $q=0$. As above, fix a basis
$X = \{x_1,\dots,x_d \}$ of $V$, ordered by $x_1 > x_2 > \dots >
x_d$, and consider the basis of $V^{\otimes n}$ consisting of the
monomials $x_{\bi} = x_{i_1}\otimes x_{i_2}\otimes\dots\otimes
x_{i_n}$ for $\bi = (i_1,i_2,\dots,i_n) \in \{1,2,\dots,d\}^n$ with
the lexicographical ordering. By equation \eqref{E:HeckDJ}, the
action of the generators $T_j$ of the Hecke algebra $\He =
\He_{n,q^2}$ on this basis is given by
\begin{equation} \label{E:HeckeKoszul}
T_j(x_{\bi}) = \begin{cases} q^2 x_{\bi} \quad&\text{if $i_j =
i_{j+1}$}\\
(q^2 - 1) x_\bi + q x_{\si_j(\bi)} &\text{if $i_j <
i_{j+1}$}\\
q x_{\si_j(\bi)} &\text{if $i_j > i_{j+1}$}
\end{cases}
\end{equation}
Here, $\si_j = (j,j+1) \in \Sy_n$ and $\si(\bi) =
(i_{\si^{-1}(1)},i_{\si^{-1}(2)},\dots,i_{\si^{-1}(n)})$ for $\si
\in \Sy_n$, as in Example~\ref{EX:SN}.

We claim that the $\He$-submodule of $V^{\otimes n}$ that is
generated by $x_\bi$ is given by
\begin{equation} \label{E:HeckeKoszul2}
\He(x_\bi) = \bigoplus_{\bi' \in \Sy_n(\bi)} \k x_{\bi'}
\end{equation}
where $\Sy_n(\bi)$ is the $\Sy_n$-orbit of  $\bi$. Indeed,
\eqref{E:HeckeKoszul} implies that each $T_\si(x_\bi)$ with $\si \in
\Sy_n$ is a linear combination of basis vectors $x_{\bi'}$ with
$\bi' \in \Sy_n(\bi)$. Hence, $\subseteq$ certainly holds in
\eqref{E:HeckeKoszul2}. For the reverse inclusion, let $\bi^*$
denote the unique non-decreasing sequence in $\Sy_n(\bi)$; so
$x_{\bi^*} = \max \{ x_{\bi'} \mid \bi' \in \Sy_n(\bi) \}$. The last
formula in \eqref{E:HeckeKoszul} implies that
\begin{equation} \label{E:HeckeKoszul3}
T(x_\bi) = q^{r(\bi)} x_{\bi^*}
\end{equation}
where $T$ is a suitable finite product of length $r(\bi) \ge 0$ in
the generators $T_j$. Since $T$ is a unit in $\He$, the inclusion
$\supseteq$ holds in \eqref{E:HeckeKoszul2}, thereby proving the
asserted equality.

Furthermore, \eqref{E:HeckeKoszul3} and \eqref{E:Heckeid3} (with $q$
replaced by $q^2$) give
\begin{equation}\label{E:ii}
q^{r(\bi)}X_n(x_\bi) = X_n(x_{\bi^*}).
\end{equation}
These elements are nonzero. For, \eqref{E:ii} implies that the
elements $X_n(x_{\bi^*})$ span the image of $X_n$ on $V^{\otimes
n}$, and their number is $\binom{d+n-1}{n}$ which is equal to the
rank of $X_n$ (cf. \cite[Eq.~(5)]{phH97}). It follows that
 $X_n(V^{\otimes n}) = \im
\rho_{n,\hR^{DJ}}(X_n)$ has a $\k$-basis consisting of the elements
\begin{equation*} 
\left\{ X_n(x_{\bi^*})  \ \big| \ \bi^* = (i_1 \le i_2 \le\dots\le
i_n) \in \{1,2,\dots,d\}^n \right\}
\end{equation*}
Next, writing
\begin{equation} \label{E:Xcoeff}
X_n(x_\bi) = \sum_{\bi' \in \Sy_n(\bi)} \lambda_{\bi'} x_{\bi'}
\end{equation}
with $\lambda_{\bi'} \in \k$, we claim that
$$
\lambda_{\si_j(\bi')} = \begin{cases} \lambda_{\bi'} \quad&\text{if $\bi' = \sigma_j(\bi')$}\\
q^{\pm 1} \lambda_{\bi'} &\text{otherwise}
\end{cases}
$$
To prove this, we may assume that $\bi' \neq \sigma_j(\bi')$. We
compute the coefficient of $x_{\sigma_j(\bi')}$ in $T_j X_n(x_\bi)$
in two ways: by \eqref{E:Heckeid3} this coefficient is equal to $q^2
\lambda_{\si_j(\bi')}$ while \eqref{E:HeckeKoszul} yields the
expression $q \lambda_{\bi'} + (q^2 - q^{1 \pm 1})
\lambda_{\si_j(\bi')}$. The claim follows from this. Writing an
arbitrary $\sigma \in \Sy_n$ as a product of the inversions
$\sigma_j$, we see that the coefficients $\lambda_{\bi'}$ in
\eqref{E:Xcoeff} only differ by a nonzero scalar, and hence they are
all nonzero since $X_n(x_\bi) \neq 0$.

By Proposition~\ref{P:confluence}, it suffices to check the
$X$-confluence conditions \eqref{E:confluence} $i=1,\ldots,N-1$. So
let $S$ be the $X$-reduction operator on $V^{\otimes N}$ with
$\Ker(S) = R$. It is easy to see from the discussion above (with
$n=N$) that $S$ is given by $S(x_{\bi^*}) =
(1-X_N/\lambda_{\bi^*})(x_{\bi^*})$ and  $S(x_\bi) = x_\bi$ for
$\bi\neq \bi^*$. According to \eqref{E:imXn2} and the discussion
above, the dimension of $(R\otimes V^{\otimes i})\cap (V^{\otimes
i}\otimes R)$ is $\binom{d+N+i-1}{N+i}$. Thus, the dimension of the
left hand side of \eqref{E:confluence} is
$d^{N+i}-\binom{d+N+i-1}{N+i}$. On the other hand  the monomials in
$V^{\otimes N+i}$ that belong to $\NRed(S_{N+i,i}) \cap
\NRed(S_{N+i,0})$ are exactly those of the form $x_{\bi^*}$ with
$\bi^* \in \{1,\ldots,d\}^{N+i}$ non-decreasing. Their number is
precisely $\binom{d+N+i-1}{N+i}$. Therefore, the dimension of
$\im(S_{N+i,i})+\im(S_{N+i,0}) = \Red(S_{N+i,i}) + \Red(S_{N+i,0})$
is at least $d^{N+i}-\binom{d+N+i-1}{N+i}$. This proves the
inequality in \eqref{E:confluence}, thereby finishing the proof of
the theorem for the case $\hR = \hR^{DJ}$.

In order to deal with an arbitrary Hecke operator $\hR$, recall that
$\He_{n,q}$ is split semisimple, having a representative set of
simple modules $M_\lambda$ indexed by the partitions $\lambda\vdash
n$; see \eqref{E:Heckess2}. We denote the representation of
$\He_{n,q}$ on $M_\lambda$ by  $\rho_\lambda$; it does not depend on
the operator $\hR$ but only on the partition $\lambda$.

Let us fix a decomposition
\begin{equation*}
V^{\otimes n}=\bigoplus_{t\in T}M_t
\end{equation*}
into simple $\He_{n,q}$-submodules $M_t$. Since all $M_t$ are
invariant under the operators $\rho_{n,\hR}(T_j)$, formula
\eqref{E:imXn2} yields the decomposition
\begin{equation*}
R_{n,i} = \bigoplus_{t\in T} \bigcap_{s=i+1}^{i+N-1}
(\rho_{n,\hR}(T_s)+1)(M_t) = \bigoplus_{t\in T} R_{n,i} \cap M_t
\end{equation*}
for all $i$. Therefore, by \cite[Lemma 1.2]{phH97}, distributivity
of the lattice $\scrT_n$ that is generated by the subspaces
$R_{n,i}$ of $V^{\otimes n}$ is equivalent to distributivity of the
lattices $\scrT_n \cap M_t \ (t \in T)$ that are generated by the
subspaces
$$
R_{n,i} \cap M_t = \bigcap_{s=i+1}^{i+N-1}
(\rho_{n,\hR}(T_s)+1)(M_t)
$$
of $M_t$.
Now, each $M_t$ is isomorphic to $M_\lambda$ for some $\lambda\vdash
n$. Therefore, the lattice $\scrT_n \cap M_t$ is isomorphic to the
lattice of subspaces of $M_\lambda$ that is generated by the
subspaces
$$
\bigcap_{s=i+1}^{i+N-1} (\rho_{\lambda}(T_s)+1)(M_\lambda)
$$
with $i+N \le n$. Finally, when $d = \dim V
> n$, then all simple $\He_{n,q}$-modules $M_\lambda$ appear in
$V^{\otimes n}$; see \cite[Proposition~5.1]{DPW91}. Thus, the
distributivity of the lattice associated to $\hR^{DJ}$, which we
have already verified, implies the distributivity of the
corresponding lattice for any Hecke operator $\hR$. This completes
the proof.
\end{proof}


\section{Koszul duality and master theorem} \label{S:duality}

In this section, $\sA = A(V,R)$ denotes an $N$-homogeneous
superalgebra that is assumed to be $N$-Koszul ($N \ge 2$).

\subsection{} \label{SS:Kdual}

By Koszulity, the complexes
\begin{equation*} 
\K(\sA)^n \colon  \ldots \to \sA_{n - \ju_N(i)} \otimes \sA^{!\
*}_{\ju_N(i)} \to \sA_{n - \ju_N(i-1)} \otimes \sA^{!\
*}_{\ju_N(i-1)} \to \dots {\to} \sA_n \to 0
\end{equation*}
are exact for $n > 0$. This yields equations in the Grothendieck
ring $R_{\eA}$ of the category $\com_{\eA}$\,:
\begin{equation} \label{E:koszul1}
\sum_{i \ge 0} (-1)^i [\sA_{n - \ju_N(i)}][ \sA^{!\ *}_{\ju_N(i)}] =
0 \qquad (n > 0)
\end{equation}
In the power series ring $R_{\eA}\llbracket t \rrbracket$ over the
Grothendieck ring R$_{\eA}$, define the Poincar{\'e} series
\begin{equation*}
P_{\sA}(t) = \sum_{n \ge 0} [\sA_{n}] t^n \qquad \text{and} \qquad
P_{\sA^{\!!*}}(t) = \sum_{n \ge 0} [A^{!\ *}_{n}] t^n
\end{equation*}
For any power series $P(t) = \sum_{n} a_n t^n$, we use the notation
\begin{equation*}
P_N(t):= \sum_{n\equiv 0,1\bmod N} (-1)^{\alpha_N(n)} a_n t^n
\end{equation*}
where $\alpha_N(n) = n - (n \bmod N)$ denotes the largest multiple
of $N$ less than or equal to $n$. Thus, $P_2(t) = P(t)$ and in
general
\begin{equation}  \label{PN(-t)}
\begin{aligned}
P_N(-t) &= \sum_{n\equiv 0,1\bmod N} (-1)^{n \bmod N} a_n t^n \\
&= \sum_{i \ge 0} (-1)^i a_{\ju_N(i)} t^{\ju_N(i)}
\end{aligned}
\end{equation}
In particular,
\begin{equation*} \label{E:PA!*}
P_{\sA^{\!!*},N}(-t) =  \sum_{i \ge 0} (-1)^i [\sA^{!\
*}_{\ju_N(i)}] t^{\ju_N(i)}
\end{equation*}
Equations \eqref{E:koszul1} are equivalent to the following Koszul
duality formula:

\begin{prop} \label{P:koszul}
For any $N$-homogeneous Koszul superalgebra $\sA$, the identity
$$
P_{\sA}(t) P_{\sA^{\!!*},N}(-t) = 1
$$
holds in $R_{\eA}\llbracket t \rrbracket$.
\end{prop}

Applying the ring homomorphism $\schi\llbracket t \rrbracket \colon
R_{\eA}\llbracket t \rrbracket \to (\eA)_{\bar 0}\llbracket t
\rrbracket$, where $\schi$ is the supercharacter map as in
Corollary~\ref{C:grothendieck}, the formula in
Proposition~\ref{P:koszul} takes the following form in $(\eA)_{\bar
0}\llbracket t \rrbracket$:

\begin{cor}
\label{C:koszul}
$\displaystyle \left( \sum_{\ell} \schi_{\sA_\ell} t^\ell \right)
\cdot \left( \sum_{m \equiv 0,1 \bmod N} (-1)^{m \bmod N}
\schi_{\sA^{!\,*}_{m}} t^{m} \right) = 1$
\end{cor}

Analogous formulas hold with the supercharacter $\schi$ replaced by
the ordinary character $\chi$ or by one of the dimensions $\dim$ and
$\sdim$.

By \eqref{E:delta3} the coaction of $\eA$ on $\sA$ sends $\sA_n$ to
$\sA_n \otimes (\eA)_n$. A similar remark holds for the
$\eA$-coaction on $\sA^{!\,*}$; see \S\ref{SS:gradeddual}.
Therefore, both factors in Corollary~\ref{C:koszul} actually belong
to the Rees subring $\prod_{n\ge 0} B_n t^n$ of $B\llbracket t
\rrbracket$, where we have put $B = (\eA)_{\bar 0}$.

\begin{example}
As an application of the Hilbert series version of
Corollary~\ref{C:koszul}, we see that the duals $\sA^!$ of the
Yang-Mills algebras $\sA = \YM^{p|q}$ are never $3$-Koszul. In fact,
by Example~\ref{EX:YM!}, we have $H_{\sA^!}(t) =
1+dt+d^2t^2+dt^3+t^4$ if $p=0$ or $q=0$ and $H_{\sA^!}(t) =
1+dt+d^2t^2+dt^3$ otherwise. In either case, $H_{\sA^!}(t)^{-1}$ has
a nonzero coefficient at $t^5$, which rules out Koszulity.
\end{example}

\subsection{A master theorem modeled on the $N$-symmetric superalgebra $\Sym_N(V)$}
\label{SS:MTSN}

We put $\sA = \Sym_N(V)$ and use the notation of
Examples~\ref{EX:SN} and \ref{EX:SN!}. In particular, we assume that
$\ch \k = 0$ and work with a fixed basis $x_1,\dots,x_d$ of $V =
\sA_1$ so that $\dg{i} = \bar 0$ for $i \le p$ and $\dg{i} = \bar 1$
for $i > p$.

From Example~\ref{EX:SN} (see also
Proposition~\ref{P:confluence}(a)), we know that a basis of
$\sA_\ell$ is given by  the monomials $x_\bi = x_{i_1} x_{i_2}\dots
x_{i_\ell}$ for sequences $\bi = (i_1,\dots,i_\ell) \in
\{1,\dots,d\}^\ell$ such that $\bi$ has no connected subsequence
$\bj = (j_1,\dots,j_N)$ of length $N$ satisfying
$$
1 \le j_1 < \ldots < j_m \le p < j_{m+1} \le \ldots \le j_N \le d =
p+q
$$
for some $m$. Adapting notation of Etingof and Pak \cite{EP06} to
our setting, we denote this set of sequences $\bi$ by
\begin{equation} \label{E:MTSN0}
\Lambda(p|q,N)_\ell
\end{equation}
For example, $\Lambda(p|q,2)_\ell$ consists of all weakly decreasing
sequences $\bi = (i_1,\dots,i_\ell)$ with entries from
$\{1,\dots,d\}$ and such that no repetition occurs in the range
$\{p+1,\dots,d\}$.

In order to evaluate the character $\schi_{\sA_\ell}$ in
Corollary~\ref{C:koszul}, recall from \eqref{E:eAformulas} that the
coaction $\delta_{\sA} \colon \sA \to \sA \otimes \eA$ is given on
the generators $x_i$ of $\sA$ by
$$
\delta_{\sA}(x_i) = \sum_j x_j \otimes z^j_i \in \sA \otimes \eA
$$
where $z^j_i = x^j \otimes x_i$ are the canonical generators of the
algebra $\eA$. For $\bi = (i_1,\dots,i_\ell) \in
\Lambda(p|q,N)_\ell$, we have
$$
\delta_{\sA}(x_\bi) = \delta_{\sA}(x_{i_1}) \delta_{\sA}(x_{i_2})
\dots \delta_{\sA}(x_{i_\ell}) \in \sA_\ell \otimes \eA
$$
Since $\sA_\ell \otimes \eA = \bigoplus_{\bi \in
\Lambda(p|q,N)_\ell} x_\bi \otimes \eA$, we can define $Z(\bi) \in
(\eA)_{\bar 0}$ by
$$
\delta_{\sA}(x_\bi) = x_\bi \otimes Z(\bi) + \left(\text{terms
supported on $\Lambda(p|q,N)_\ell \setminus \{\bi\}$}\right)
$$
Then \eqref{E:sch4} becomes
\begin{equation} \label{E:MTSN1}
\schi_{\sA_\ell} = \sum_{\bi \in \Lambda(p|q,N)_\ell}
(-1)^{\dg{\bi}} Z(\bi)
\end{equation}
with $\dg{\bi} = \dg{i_1} + \dots + \dg{i_\ell}$.

Now consider the super bialgebra $\sB = \sO(\Endo(V)) = \k[x^i_j
\mid 1 \le i,j \le d]$ defined in \S\ref{SSS:EV2} and recall that
the $x^i_j$ are supercommuting variables of parity $\dg{i} + \dg{j}$
over $\k$. Restricting the comodule $\sA_\ell$ to $\sB$ along the
map $\varphi \colon \en \Sym_N(V) \to \sB$, $z^i_j \mapsto x^i_j$ in
\eqref{E:endSN} we must replace $Z(\bi)$ in \eqref{E:MTSN1} by
$X(\bi):= \varphi(Z(\bi)) \in \sB_{\bar 0}$. Thus, writing
\begin{equation*}
y_i = \sum_j x_j \otimes x^j_i \in \sA \otimes \sB
\end{equation*}
and $y_\bi = y_ {i_1}\dots y_{i_\ell} \in \sA_\ell \otimes \sB =
\bigoplus_{\bj \in \Lambda(p|q,N)_\ell} x_\bj \otimes \sB$ for $\bi
= (i_1,\dots,i_\ell)$, we have
\begin{equation} \label{E:MTSN2}
y_\bi = x_\bi \otimes X(\bi) + \left(\text{terms supported on
$\Lambda(p|q,N)_\ell \setminus \{\bi\}$}\right)
\end{equation}

As for the supercharacter of $\sA^{!\,*}_{m}$, recall from
\eqref{E:A!*} and \eqref{E:SN!} that, for all $n \ge N$,
\begin{equation*}
\sA^{!\,*}_{n} = \bigcap_{i+j+N=n} V^{\otimes i} \otimes R \otimes
V^{\otimes j} = \Gr^n V
\end{equation*}
Viewing $\sA^{!\,*}_{n} = \Gr^n V$ as a comodule over $\sB =
\sO(\Endo(V))$, the supercharacter of $\sA^{!\,*}_{n}$ is the
$n^\text{th}$ elementary supersymmetric function $e_n$ which we
know, by Proposition~\ref{P:charfctn}, to be identical to the
coefficient at $t^n$ of the characteristic function $\ber(1+ tX)$ of
the generic supermatrix $X = \left( x^i_j \right)_{1 \le i,j \le d}$
of type $p|q$; so the diagonal blocks $X_{11} = \left( x^i_j
\right)_{1 \le i,j \le p}$ and $X_{22} = \left( x^i_j \right)_{p+1
\le i,j \le p+q}$ consist of even entries while all other entries
are odd.

To summarize, we obtain the following super-version of \cite[Theorem
2]{EP06}.

\begin{thm} \label{T:MTSN}
Let $X = \left( x^i_j \right)_{d \times d}$ be the generic
supermatrix of type $p|q$. Then
$$
\left( \sum_{\ell} \sum_{\bi \in \Lambda(p|q,N)_\ell}
(-1)^{\dg{\bi}} X(\bi)\, t^\ell \right) \cdot \left( \sum_{m \equiv
0,1 \bmod N} (-1)^{m \bmod N} e_{m} t^{m} \right) = 1
$$
holds in the power series ring $\k[x^i_j \mid \text{\rm all } i,j\,
]_{\bar 0}\llbracket t \rrbracket$. Here $\Lambda(p|q,N)_\ell$ and
$X(\bi)$ are defined by \eqref{E:MTSN0} and \eqref{E:MTSN2},
respectively, and the $e_m$ are the coefficients of the
characteristic function $\ber(1 + tX) = \sum_{n \ge 0} e_n t^n$ of
$X$.
\end{thm}

\subsection{}
\label{SS:sdimSN}

As an application of Theorem~\ref{T:MTSN} , we determine the
superdimension Hilbert series
\begin{equation*} \label{E:sdimSN1}
H^s_{\sA}(t) = \sum_{\ell \ge 0} \sdim_\k \sA_\ell\, t^\ell
\end{equation*}
for the $N$-symmetric superalgebra $\sA = \Sym_N(V)$. For
the pure even case, this was already done by Etingof and Pak
\cite{EP06} . The notations of
\S\ref{SS:MTSN} remain in effect.

In view of Corollary~\ref{C:grothendieck}, the superdimension Poincar\'e series
follows by applying the counit $\e \colon \sB \to
\k$ to the equation in Theorem~\ref{T:MTSN}. Indeed, by
\eqref{E:Bformulas}, the counit $\e$ sends $X \mapsto 1_{d \times
d}$, and hence the elements $X(\bi)$ in \eqref{E:MTSN2} all map to
$1$. Therefore, the first factor in Theorem~\ref{T:MTSN} becomes
\begin{equation*} \label{E:sdimSN2}
H^s_{\sA}(t) = \sum_{\ell \ge 0} \left( \sum_{\bi \in
\Lambda(p|q,N)_\ell} (-1)^{\dg{\bi}}\right)\,  t^\ell
\end{equation*}
For the second factor, note that
\begin{equation*}
\ber(1 + t\, 1_{d \times d}) = (1+t)^{p-q}
\end{equation*}
by \eqref{E:ber}. Thus,
\begin{equation} \label{E:sdimSN3}
\begin{aligned}
H^s_{\sA}(t) &= \sum_{\ell \ge 0} \left( \sum_{\bi
\in \Lambda(p|q,N)_\ell} (-1)^{\dg{\bi}}\right)\,  t^\ell \\
&=
\begin{cases}
\left( \displaystyle  \sum_{m \equiv 0,1 \bmod N} (-1)^{m \bmod N}
\binom{p-q}{m}\,t^m \right)^{-1}  &\text{if $p \ge q$} \\
\left( \displaystyle   \sum_{m \equiv 0,1 \bmod N} (-1)^{\alpha_N(m)}
\binom{m + q-p-1}{q-p-1}\,t^m \right)^{-1}   &\text{if $p < q$}
 \end{cases}
\end{aligned}
\end{equation}
where $\alpha_N(m) = m - (m \bmod N)$ denotes the largest multiple
of $N$ less than or equal to $m$ as in \S\ref{SS:Kdual}.

\subsection{}
\label{SS:dimSN}

The ordinary Hilbert series
$H_{\sA}(t) = \sum_{\ell \ge 0} \dim_\k \sA_\ell\, t^\ell$
of the $N$-symmetric superalgebra $\sA = \Sym_N(V)$ is as follows.
Recall from
\S\ref{SS:MTSN} that
\begin{equation*}
\dim_\k \sA_\ell = |\Lambda(p|q,N)_\ell|
\end{equation*}
and from \eqref{E:SN!dim} that
\begin{equation*} \label{}
\dim_\k \sA^!_n = \begin{cases} d^n & \text{if $n < N$}
\\ \sum_{r+s = n} \binom{p}{r} \binom{q+s-1}{s} & \text{if $n \ge N$}
\end{cases}
\end{equation*}
Therefore,
the Hilbert series is
\begin{equation} \label{E:HSN}
\begin{aligned}
H_{\sA}(t) &= \sum_{\ell \ge 0} |\Lambda(p|q,N)_\ell|\,  t^\ell \\
&= \left( \displaystyle  \sum_{m \equiv 0,1 \bmod N} (-1)^{m \bmod
N} \left( \sum_{r+s = m} \binom{p}{r} \binom{q+s-1}{s} \right)\,t^m
\right)^{-1}
\end{aligned}
\end{equation}

\subsection{}
\label{SS:dimSRN}

Less is known about the Hilbert series of the $N$-homogeneous
superalgebras $\sA = \Gr_{\hR,N} $ associated to an arbitrary Hecke
operator $\hR \colon V^{\otimes 2} \to V^{\otimes 2}$ on a vector
superspace $V$ ; see Example~\ref{EX:SRN}. Recall that $\sA =
A(V,R)$ with $R = \im \rho_{\hR}(X_N) \subseteq V^{\otimes N}$. For
any $N$-homogeneous algebra $\sA = A(V,R)$, we have
\begin{equation*}
\dim_\k \sA^!_n = \dim_\k \bigcap_{i+j+N=n}  V^{\otimes j} \otimes R
\otimes V^{\otimes i}
\end{equation*}
by \eqref{E:Rn!} and \eqref{E:Rn!2}. For $R = \im \rho_{\hR}(X_N)$
in particular, \eqref{E:imXn} further implies that
\begin{equation*}
\bigcap_{i+j+N=n}  V^{\otimes j} \otimes R \otimes V^{\otimes i} =
\rho_{\hR}(X_n)\left(V^{\otimes n}\right)
\end{equation*}
holds for $n \ge N$. Now \cite[Theorem 3.5]{phH99} implies that
\begin{equation*}
H_{\Gr_{\hR,2}{}^!}(t) = \frac{\prod_{\ell=1}^r (1+a_\ell t)}{\prod_{m=1}^s
(1-b_m t)}
\end{equation*}
where $(r,s)$ is the birank of $\hR$ and  $a_\ell$ and $b_m$ are positive real numbers.
For example, in the situation of \ref{SS:dimSN}, $(r,s)=(p,q)$ and $a_\ell=b_m=1$.

 For any
complex power series $P(t)$, the power series $P_N(-t)$ in
\eqref{PN(-t)} can be written as
\begin{equation*}
P_N(-t) = \frac{1}{N} \sum_{i=1}^{N-1} (1 - \zeta_N^{-i})P(\zeta_N^i
t)
\end{equation*}
where $\zeta_N = e^{2\pi i/N}$. In particular,
\begin{equation*}
\begin{aligned}
H_{\sA^{\!!*},N}(-t)  &= \frac{1}{N} \sum_{i=1}^{N-1} (1 -
\zeta_N^{-i}) \frac{\prod_{\ell=1}^r (1+a_\ell \zeta_N^i
t)}{\prod_{m=1}^s (1-b_m \zeta_N^i t)} \\
&= \frac{Q_{N,\bf a,b}(t)}{\prod_{m=1}^s (1 + b_mt + \ldots + b_m^{N-1}t^{N-1})}
\end{aligned}
\end{equation*}
for some real polynomial $Q_{N,\bf a,\bf b}(t)$ with coefficients being polynomial in ${\bf a}=(a_\ell)$ and ${\bf b}=(b_m)$. Therefore, the Hilbert series of
$\sA$ has the form
\begin{equation}
H_{\sA}(t) = \frac{\prod_{m=1}^s (1 + b_mt + \ldots +
b_m^{N-1}t^{N-1})}{Q_{N,\bf a,b}(t)}
\end{equation}
Notice that the fraction on the right-hand side is reduced.

In particular, \eqref{E:HSN} has the form
\begin{equation}
H_{\sA}(t)=\frac{(1-t^N)^s}{(1-t)^sQ_{N,\bf 1,\bf 1}(t)}
\end{equation}


\section*{Appendix}

\setcounter{section}{1}
\renewcommand\thesection{\Alph{section}}%

For lack of a suitable reference, we include here a proof of
Proposition~\ref{P:intro} that was stated in the Introduction. Our
proof is based on the proof of \cite[Proposition 2.1]{rBnM06} and on
additional details that were communicated to us by Roland Berger.
For the basics concerning graded algebras, we refer the reader to
\cite[Chap.~II \S11]{nB70} or \cite{rB05}.

As in the Introduction, $\sA = \bigoplus_{n \ge 0} \sA_n$ denotes an
arbitrary connected $\Z_{\ge 0}$-graded $\k$-algebra and $V$ is a
graded subspace of $\sA_+ = \bigoplus_{n > 0} \sA_n$ satisfying
$\sA_+ = V \oplus \sA_+^2$. Thus, $\T(V)/I \iso \sA$ for some graded
ideal $I$ of $\T(V)$. For convenience, we state
Proposition~\ref{P:intro} again:

\begin{prop*}
The relation ideal $I$ of $\sA$ lives in degrees $\ge N$ if and only
if $\Tor^{\sA}_i(\k,\k)$ lives in degrees $\ge \ju_N(i) =
\begin{cases} \tfrac{i}{2}N \quad &\text{if $i$ is even}\\
\tfrac{i-1}{2}N + 1 \quad &\text{if $i$ is odd}
\end{cases}$
\end{prop*}

\begin{proof} Let
\begin{equation*}
P \colon \quad \dots \to P_i \stackrel{d_i}{\tto} P_{i-1}
\stackrel{d_{i-1}}{\tto} \dots \stackrel{d_1}{\tto} P_0
\stackrel{d_0}{\tto} \k \to 0
\end{equation*}
be a minimal graded-free resolution of the trivial left $\sA$-module
$\k$. Thus, all $P_i$ have the form $P_i = \sA \otimes E_i$ for some
graded subspace $E_i \subseteq \Ker d_{i-1}$ which is chosen so that
\begin{equation} \label{E:E}
\Ker d_{i-1} = E_i \oplus \sA_+\Ker d_{i-1}
\end{equation}
In particular, we may take $E_0 = \k$ and $E_1 = V$. The
differential $d_i \colon P_i \to P_{i-1}$ is the graded $\sA$-module
map that is defined by the inclusion $E_i \into P_{i-1}$. By the
graded Nakayama Lemma (e.g., \cite[p.~AII.171, Prop.~6]{nB70}), our
choice of $E_i$ implies that
\begin{equation} \label{E:exact}
\im d_{i} = \sA E_i = \Ker d_{i-1}\quad \text{and}\quad \Ker d_i
\subseteq \sA_+ \otimes E_i = \sA_+ P_i
\end{equation}
for all $i$. Consequently, the complex $\k \otimes_{\sA} P$ has zero
differential, and hence
\begin{equation*}
\Tor^{\sA}_i(\k,\k) \cong \k \otimes_{\sA} P_i \cong E_i
\end{equation*}
In particular,
\begin{equation*} \label{E:0and1}
\Tor^{\sA}_0(\k,\k) \cong \k \quad \text{and} \quad
\Tor^{\sA}_1(\k,\k) \cong V = \sA_+/\sA_+^2
\end{equation*}
live in degrees $0 = \ju_N(0)$ and $\ge 1 = \ju_N(1)$, respectively.
Moreover, the kernel of $d_1 \colon P_1 = \left(\T(V)/I\right)
\otimes V \to P_0 = \sA$ is exactly $I/I \otimes V$, and so
\begin{equation*} \label{E:2}
\Tor^{\sA}_2(\k,\k) \cong \Ker d_1/\sA_+\Ker d_1 \cong I/\left(V
\otimes I + I \otimes V\right)
\end{equation*}
Therefore, $I$ lives in degrees $\ge N$ if and only if
$\Tor^{\sA}_2(\k,\k)$ lives in degrees $\ge N = \ju_N(2)$.

For the remainder of the proof, assume that  $I$ lives in degrees
$\ge N$. We will show by induction on $i$ that $\Tor^{\sA}_i(\k,\k)
= E_i$ lives in degrees $\ge \ju_N(i)$ for all $i$. The cases $i \le
2$ have been checked above. Assume that $E_i$ lives in degrees $\ge
\ju_N(i)$ and similarly for $E_{i-1}$. By \eqref{E:exact}, we know
that $E_{i+1} \subseteq \Ker d_i \subseteq \sA_+ \otimes E_{i}$ and
so $E_{i+1}$ certainly lives in degrees $\ge \ju_N(i) +1$. Since
$\ju_N(i) +1 = \ju_N(i +1)$ when $i$ is even (or when $i$ is
arbitrary and $N=2$), we are done in these cases. From now on, we
assume that $i$ is odd. We must show that $E_{i+1}$ lives in degrees
$\ge \ju_N(i+1) = \tfrac{i+1}{2}N$. Since $E_{i+1} \subseteq \Ker
d_i$, it suffices to show that $d_i$ is injective in degrees $<
\tfrac{i+1}{2}N$, and since $E_i$ lives in degrees $\ge \ju_N(i)=
\tfrac{i-1}{2}N+1$, our goal is to show that $d_i$ is injective on
all homogeneous components $P_{i,n}$ of $P_i$ in degrees $n =
\tfrac{i-1}{2}N+j$ with $j = 1,\dots,N-1$. Put $m = \tfrac{i-1}{2}N$
for simplicity and note that
\begin{equation} \label{E:i}
P_{i,m+j} = \bigoplus_{\ell = 1}^j \sA_{j-\ell}\otimes E_{i,m+\ell}
\end{equation}
and
\begin{equation} \label{E:i-1}
P_{i-1,m+j} = \bigoplus_{k = 0}^j \sA_{j-k}\otimes E_{i-1,m+k}
\end{equation}
since $E_{i-1}$ lives in degrees $\ge \ju_N(i-1) = m$. The
proposition will be a consequence of the following claims:
\begin{enumerate}
\item $d_i$ is injective
on all summands $\sA_{j-\ell}\otimes E_{i,m+\ell}$ in \eqref{E:i},
and
\item the subspaces $d_i\left(\sA_{j-\ell}\otimes E_{i,m+\ell}\right)
= \sA_{j-\ell}E_{i,m+\ell}$ for $\ell = 1,\dots,j$ form a direct sum
inside $P_{i-1,m+j}$.
\end{enumerate}

In order to prove (a), recall that the restriction of $d_i$ to
$E_{i,m+\ell}$ is the inclusion
\begin{equation*}
E_{i,m+\ell} \into
P_{i-1,m+\ell} = \bigoplus_{k = 0}^\ell \sA_{\ell-k}\otimes
E_{i-1,m+k}
\end{equation*}
Hence, the effect of $d_i$ on the $\ell^\text{th}$ summand in
\eqref{E:i} is the embedding
\begin{equation*}
\sA_{j-\ell}\otimes E_{i,m+\ell} \into \bigoplus_{k = 0}^\ell
\sA_{j-\ell} \otimes \sA_{\ell-k}\otimes E_{i-1,m+k}
\end{equation*}
followed by the map
\begin{equation*}
\bigoplus_{k = 0}^\ell \sA_{j-\ell} \otimes \sA_{\ell-k}\otimes
E_{i-1,m+k} \tto \bigoplus_{k = 0}^\ell \sA_{j-k}\otimes E_{i-1,m+k}
\subseteq P_{i-1,m+j}
\end{equation*}
which is given by the multiplication map $\sA_{j-\ell} \otimes
\sA_{\ell-k} \to \sA_{j-k}$. Since $j-k < N$, our hypothesis on $I$
implies that $\sA_{j-k} \cong \T(V)_{j-k}$, and similarly
$\sA_{j-\ell} \cong \T(V)_{j-\ell}$ and $\sA_{\ell-k} \cong
\T(V)_{\ell-k}$. Therefore, the above multiplication map is
identical with the injection $\T(V)_{j-\ell} \otimes \T(V)_{\ell-k}
\into \T(V)_{j-k}$ in $\T(V)$. This proves (a).

For (b), we proceed by induction on $j$. The case $j=1$ being
obvious, let $1 \le j \le N-2$ and assume that (ii) holds for
$1,\dots,j$. We wish to show that the subspaces
$\sA_{j+1-\ell}E_{i,m+\ell}$ $(\ell = 1,\dots,j+1)$ of
$P_{i-1,m+j+1}$ form a direct sum. First, by \eqref{E:E} we have
$E_{i,m+j+1} \cap \sA_+ \Ker d_{i-1} = 0$ while $\sum_{\ell = 1}^j
\sA_{j+1-\ell}E_{i,m+\ell} \subseteq \sA_+ \Ker d_{i-1}$. Therefore,
it suffices to show that the sum $\sum_{\ell = 1}^j
\sA_{j+1-\ell}E_{i,m+\ell}$ is direct. To this end, note that
$\sA_{j+1-\ell} = \sum_{d \ge 1} V_d\sA_{j+1-d-\ell}$ holds for all
$\ell \le j$. Hence,
\begin{equation*}
\sum_{\ell = 1}^j \sA_{j+1-\ell}E_{i,m+\ell} = \sum_{d \ge 1} V_d
\sum_{\ell = 1}^j \sA_{j+1-d-\ell}E_{i,m+\ell}
\end{equation*}
By induction, $\sum_{\ell = 1}^j \sA_{j+1-d-\ell}E_{i,m+\ell}$ is a
direct sum inside $P_{i-1,m+j+1-d}$. Thus, it suffices to show that
the sum $\sum_{d \ge 1} V_d P_{i-1,m+j+1-d} \subseteq P_{i-1,m+j+1}$
is direct. But \eqref{E:i-1} gives
\begin{equation*}
P_{i-1,m+j+1} = \bigoplus_{k = 0}^{j+1} \sA_{j+1-k}\otimes
E_{i-1,m+k} = \bigoplus_{k = 0}^{j+1} \T(V)_{j+1-k}\otimes
E_{i-1,m+k}
\end{equation*}
where the last equality holds since all $j+1-k < N$. Therefore,
\begin{equation*}
\sum_{d \ge 1} V_d P_{i-1,m+j+1-d} = \bigoplus_{d \ge 1} V_d \otimes
\bigoplus_{k = 0}^{j+1-d} \T(V)_{j+1-d-k}\otimes E_{i-1,m+k}
\end{equation*}
as desired. This proves (b), thereby completing the proof of the
proposition.
\end{proof}


\begin{ack}
The authors wish to thank Roland Berger for his helpful comments
throughout the completion of this paper. Work on this article was
initiated during a visit of ML to the Universit{\'e} Montpellier 2
in June 2006. ML would like to thank Claude Cibils for arranging
this visit and for his warm hospitality in Montpellier, and James
Zhang for help with some examples in this article.
\end{ack}


\def\cprime{$'$}
\providecommand{\bysame}{\leavevmode\hbox
to3em{\hrulefill}\thinspace}
\providecommand{\MR}{\relax\ifhmode\unskip\space\fi MR }
\providecommand{\MRhref}[2]{%
  \href{http://www.ams.org/mathscinet-getitem?mr=#1}{#2}
} \providecommand{\href}[2]{#2}



\begin{thebibliography}{10}

\bibitem{mAwS87}
Michael Artin and William~F. Schelter, \emph{Graded algebras of
global dimension {$3$}}, Adv. in Math. \textbf{66} (1987), no.~2,
171--216. \MR{MR917738 (88k:16003)}

\bibitem{BGS96}
Alexander Beilinson, Victor Ginzburg, and Wolfgang Soergel,
\emph{Koszul duality patterns in representation theory}, J. Amer.
Math. Soc. \textbf{9} (1996), no.~2, 473--527. \MR{MR1322847
(96k:17010)}

\bibitem{fB87}
Felix~Alexandrovich Berezin, \emph{Introduction to superanalysis},
Mathematical Physics and Applied Mathematics, vol.~9, D. Reidel
Publishing Co., Dordrecht, 1987, Edited and with a foreword by A. A.
Kirillov, With an appendix by V. I. Ogievetsky, Translated from the
Russian by J. Niederle and R. Koteck\'y, Translation edited by
Dimitri Le\u\i tes. 

\bibitem{rB98}
Roland Berger, \emph{Confluence and {K}oszulity}, J. Algebra
\textbf{201} (1998), no.~1, 243--283. 

\bibitem{rB01}
\bysame, \emph{Koszulity for nonquadratic algebras}, J. Algebra
\textbf{239} (2001), no.~2, 705--734. 

\bibitem{rB05}
\bysame, \emph{Dimension de {H}ochschild des alg\`ebres gradu\'ees},
C. R. Math. Acad. Sci. Paris \textbf{341} (2005), no.~10, 597--600.

\bibitem{rBmDVmW}
Roland Berger, Michel Dubois-Violette, and Marc Wambst,
\emph{Homogeneous algebras}, J. Algebra \textbf{261} (2003), no.~1,
172--185. 

\bibitem{rBnM06}
Roland Berger and Nicolas Marconnet, \emph{Koszul and {G}orenstein
properties for homogeneous algebras}, Algebr. Represent. Theory
\textbf{9} (2006), no.~1, 67--97. 

\bibitem{nB70}
Nicolas Bourbaki, \emph{{A}lg\`ebre, {C}hapitres 1 \`a 3}, Hermann,
Paris, 1970. 

\bibitem{aCmDV02}
Alain Connes and Michel Dubois-Violette, \emph{Yang-{M}ills
algebra}, Lett. Math. Phys. \textbf{61} (2002), no.~2, 149--158.

\bibitem{aCmDV03}
\bysame, \emph{{Yang-Mills and some related algebras}}, 2004,
arXiv:math-ph/0411062.

\bibitem{pDjM}
Pierre Deligne and James~S. Milne, \emph{Tannakian categories},
Hodge cycles, motives, and Shimura varieties, Lecture Notes in
Math., vol. 900, Springer-Verlag, Berlin-New York, 1982,
pp.~101--228. 

\bibitem{rDgJ86}
Richard Dipper and Gordon James, \emph{Representations of {H}ecke
algebras of general linear groups}, Proc. London Math. Soc. (3)
\textbf{52} (1986), no.~1, 20--52. 

\bibitem{rDgJ91}
\bysame, \emph{{$q$}-tensor space and {$q$}-{W}eyl modules}, Trans.
Amer. Math. Soc. \textbf{327} (1991), no.~1, 251--282.

\bibitem{DPW91}
Jie Du, Brian Parshall, and Jian~Pan Wang, \emph{Two-parameter
quantum linear groups and the hyperbolic invariance of {$q$}-{S}chur
algebras}, J. London Math. Soc. (2) \textbf{44} (1991), no.~3,
420--436. 

\bibitem{EP06}
Pavel Etingof and Igor Pak, \emph{{An algebraic extension of the
MacMahon Master Theorem}}, Proc. Amer. Math. Soc. (to appear),
arXiv:math.CO/0608005.

\bibitem{FH06a}
Dominique Foata and Guo-Niu Han, \emph{{A basis for the right
quantum algebra and the "1=q" principle}}, arXiv:math.CO/0603463.

\bibitem{FH06b}
\bysame, \emph{{A New Proof of the Garoufalidis-Le-Zeilberger
Quantum MacMahon Master Theorem}}, arXiv:math.CO/0603464.

\bibitem{FH06c}
\bysame, \emph{{Specializations and Extensions of the quantum
MacMahon Master Theorem}}, arXiv:math.CO/0603466.

\bibitem{GLZxx}
Stavros Garoufalidis, Thang TQ~Le, and Doron Zeilberger, \emph{The
quantum {M}ac{M}ahon {M}aster {T}heorem}, Proc. Natl. Acad. of Sci.
\textbf{103} (2006), 13928--13931, arXiv:math.QA/0303319.

\bibitem{eGetal04}
E.~L. Green, E.~N. Marcos, R.~Mart{\'{\i}}nez-Villa, and Pu~Zhang,
\emph{{$D$}-{K}oszul algebras}, J. Pure Appl. Algebra \textbf{193}
(2004), no.~1-3, 141--162. 

\bibitem{dG90}
D.~I. Gurevich, \emph{Algebraic aspects of the quantum
{Y}ang-{B}axter equation}, Algebra i Analiz \textbf{2} (1990),
no.~4, 119--148. 

\bibitem{GPS05}
D.~I. Gurevich, P.~N. Pyatov, and P.~A. Saponov, \emph{The
{C}ayley-{H}amilton theorem for quantum matrix algebras of {${\rm
GL}(m\vert n)$} type}, Algebra i Analiz \textbf{17} (2005), no.~1,
160--182. 

\bibitem{GPS06}
\bysame, \emph{Quantum matrix algebras of {${\rm GL}(m\vert
n)$}-type: the structure of the characteristic subalgebra and its
spectral parametrization}, Teoret. Mat. Fiz. \textbf{147} (2006),
no.~1, 14--46. 

\bibitem{phH97}
Ph{\`u}ng~H{\^o} Hai, \emph{Koszul property and {P}oincar\'e series
of matrix bialgebras of type {$A\sb n$}}, J. Algebra \textbf{192}
(1997), no.~2, 734--748. 

\bibitem{phH99}
\bysame, \emph{Poincar\'e series of quantum spaces associated to
{H}ecke operators}, Acta Math. Vietnam. \textbf{24} (1999), no.~2,
235--246. 

\bibitem{phH02a}
\bysame, \emph{Realizations of quantum hom-spaces, invariant theory,
and quantum determinantal ideals}, J. Algebra \textbf{248} (2002),
no.~1, 50--84. 

\bibitem{HL}
Ph{\`u}ng~H{\^o} Hai and Martin Lorenz, \emph{{Koszul algebras and
the quantum MacMahon Master Theorem}}, Bull. London Math. Soc. (to
appear), arXiv:math.QA/0603169.


\bibitem{cK95}
Christian Kassel, \emph{Quantum groups}, Graduate Texts in
Mathematics, vol. 155, Springer-Verlag, New York, 1995.

\bibitem{hKtV05}
H.~M. Khudaverdian and Th.~Th. Voronov, \emph{Berezinians, exterior
powers and recurrent sequences}, Lett. Math. Phys. \textbf{74}
(2005), no.~2, 201--228. 

\bibitem{mK07a}
Matja{\v z} Konvalinka, \emph{{A generalization of Foata's
fundamental transformation and its applications to the right-quantum
algebra}}, preprint (March 2007), arXiv:math/0703203v1 [math.CO].

\bibitem{mK07b}
\bysame, \emph{{Non-commutative Sylvester's determinantal
identity}}, preprint (March 2007), arXiv:math/0703213v1 [math.CO].

\bibitem{mKiP}
Matja{\v z} Konvalinka and Igor Pak, \emph{Noncommutative extensions
of the {M}ac{M}ahon {M}aster {T}heorem}, Adv. in Math. (to appear);
arXiv:math.CO/0607737.

\bibitem{igMac}
I.~G. Macdonald, \emph{Symmetric functions and {H}all polynomials},
second ed., Oxford Mathematical Monographs, The Clarendon Press
Oxford University Press, New York, 1995, With contributions by A.
Zelevinsky, Oxford Science Publications. 

\bibitem{paMM60}
Percy~A. MacMahon, \emph{Combinatory analysis}, Two volumes (bound
as one), Chelsea Publishing Co., New York, 1960, reprint (1960) in
one volume of two volumes originally published by Cambridge Univ.
Press, 1915 and 1916. 

\bibitem{yM87}
Yu.~I. Manin, \emph{Some remarks on {K}oszul algebras and quantum
groups}, Ann. Inst. Fourier (Grenoble) \textbf{37} (1987), no.~4,
191--205. 

\bibitem{yM88}
\bysame, \emph{Quantum groups and noncommutative geometry},
Universit\'e de Montr\'eal Centre de Recherches Math\'ematiques,
Montreal, QC, 1988. 

\bibitem{yM89}
Yuri~I. Manin, \emph{Multiparametric quantum deformation of the
general linear supergroup}, Comm. Math. Phys. \textbf{123} (1989),
no.~1, 163--175. 

\bibitem{yM91}
\bysame, \emph{Topics in noncommutative geometry}, M. B. Porter
Lectures, Princeton University Press, Princeton, NJ, 1991.

\bibitem{yM97}
\bysame, \emph{Gauge field theory and complex geometry}, second ed.,
Grundlehren der Mathematischen Wissenschaften [Fundamental
Principles of Mathematical Sciences], vol. 289, Springer-Verlag,
Berlin, 1997, Translated from the 1984 Russian original by N.
Koblitz and J. R. King, With an appendix by Sergei Merkulov.

\bibitem{PP05}
Alexander Polishchuk and Leonid Positselski, \emph{Quadratic
algebras}, University Lecture Series, vol.~37, American Mathematical
Society, Providence, RI, 2005. 

\bibitem{dStjZ97}
Darin~R. Stephenson and James~J. Zhang, \emph{Growth of graded
{N}oetherian rings}, Proc. Amer. Math. Soc. \textbf{125} (1997),
no.~6, 1593--1605. 

\bibitem{gT04}
Gijs~M. Tuynman, \emph{Supermanifolds and supergroups}, Mathematics
and its Applications, vol. 570, Kluwer Academic Publishers,
Dordrecht, 2004. 

\bibitem{mW93}
Marc Wambst, \emph{Complexes de {K}oszul quantiques}, Ann. Inst.
Fourier (Grenoble) \textbf{43} (1993), no.~4, 1089--1156.

\end{thebibliography}
\end{document}